\documentclass[11pt]{article}

\usepackage{epsfig}
\usepackage{amssymb,amsmath}

\usepackage{latexsym}

\renewcommand{\H}{{\bf H}}
\newcommand{\R}{{\mathbb R}}
\newcommand{\pp}{{\bf p}}

\newcommand{\V}{{\bf V}} 
\newcommand{\I}{{\bf I}}
\newcommand{\lr}{\leftrightarrow}
\newtheorem{theorem}{Theorem}[section]
\newtheorem{proposition}[theorem]{ Proposition}

\newtheorem{lemma}[theorem]{Lemma}

\newtheorem{remarks}[theorem]{Remarks}

\begin{document}
\title{Vertices and inflexions of plane sections of surfaces in ${\mathbb R}^3$}
\author{Andr\'{e} Diatta  and Peter Giblin \\
University of Liverpool, Liverpool L69 3BX, England \\ email {\tt adiatta@liv.ac.uk  pjgiblin@liv.ac.uk}}
\date{}
\maketitle             

\begin{abstract}
We discuss the behaviour of vertices and inflexions of  one-parameter families of plane curves which include a singular member. 
These arise a sections of smooth surfaces by families of planes parallel to the tangent plane at a given point.
We cover all the generic cases, namely elliptic ($A_1$), umbilic ($A_1$), hyperbolic ($A_1$), parabolic ($A_2$) and cusp of Gauss ($A_3$) points.
This work is preliminary to an investigation of symmetry sets and medial axes for these families of curves, reported elsewhere.
 \end{abstract}

\section{Introduction}
\label{s:intro} 

Let $M$ be a smooth surface, and {\bf p} be a point of $M$. We shall consider the intersection of $M$ with a family
of planes parallel to the tangent plane at {\bf p}. This family of plane curves contains a singular member, when the
plane is the tangent plane itself; generically the other members of the family close to the tangent plane are nonsingular curves.

The motivation for this work comes from computer vision, where the surface is the intensity surface $z=f(x,y)$
corresponding to the intensity function $f$ of a two-dimensional image, and the plane curves are level
sets of this function, that is isophotes. A great deal of information about the shape of these level sets and the
way they evolve through the singular level set is contained in the family of so-called {\em symmetry sets} 
and {\em medial axes} of the
level sets (see for example \cite{MathsOfSurfaces2000}). These sets in  turn take some of their structure from the pattern of
vertices and inflexions (curvature extrema and zeros) of the level set. 

In this article we concentrate on the
vertices and inflexions, and apply this and other results to the study of symmetry sets in articles to appear
elsewhere~\cite{scalespace05,maths-of-surfaces}. Besides the patterns of vertices and inflexions we also study the limiting curvatures
at the vertices as the level set approaches the singular member of the family. 

The contact between a surface and its tangent plane at {\bf p} is an {\em affine invariant} of the surface.
Likewise the inflexions on the intersections with nearby planes are affine invariants, but we are
also interested in the curvature extrema on these sections, and these are {\em euclidean}
invariants.
For a generic surface $M$, the contact between the surface and its tangent plane at a point {\bf p},
as measured by the height function in the normal direction at {\bf p}, can be of the following
types.
See for example \cite{solid-shape} for the geometry of these situations, and \cite{b-g-t,bgt95,mumford-book}
for an extensive discussion of the singularity theory.\\
$\bullet$ \ \ The contact at {\bf p} is ordinary (`$A_1$ contact'), at 
an elliptic point or at a hyperbolic point (occupying {\em regions} of $M$). The intersection of $M$ with its tangent
plane at {\bf p} is locally an isolated point or a pair of transverse smooth arcs. (As regards 
contact there is no distinction between `ordinary' elliptic points and umbilics, where
the principal curvatures coincide. But as we shall see there is a great deal of difference when
we consider vertices of the plane sections.) \\
$\bullet$ \ \ The contact is of type $A_2$ at  parabolic points (generically forming 
{\em smooth curves} on $M$), where the asymptotic directions coincide.
 The intersection
of $M$ with its tangent plane at {\bf p} is locally a cusped curve. \\
$\bullet$ \ \ The contact is of type $A_3$ at a cusp of Gauss, where the parabolic curve is tangent to
the asymptotic direction at {\bf p} (these are {\em isolated points} of $M$).
There are two types, the elliptic cusp and the hyperbolic cusp. The intersection
with the tangent plane is locally an isolated point or a pair of tangential arcs.

\medskip

Other authors have considered the $A_1$ cases, using different techniques and with slightly different motivations
from ours. Vertices in the $A_1$ case are studied by Uribe-Vargas in~\cite{uribe-vargas} 
and inflexions in the same case by Garay in~\cite{garay},
using more sophisticated techniques of singularity theory aimed at finding normal forms up
to an appropriate equivalence.
Our very detailed results on the other hand combine vertices and inflexions and 
apply to all three cases $A_1, A_2, A_3$ above. They are obtained by direct calculation: our motivation,
as above, is to facilitate investigation of the symmetry sets of surface sections, and we do not as yet know
how to fit this into a more general theory.

\medskip

\begin{figure}
  \begin{center}
  \leavevmode
\epsfxsize=1in
\epsffile{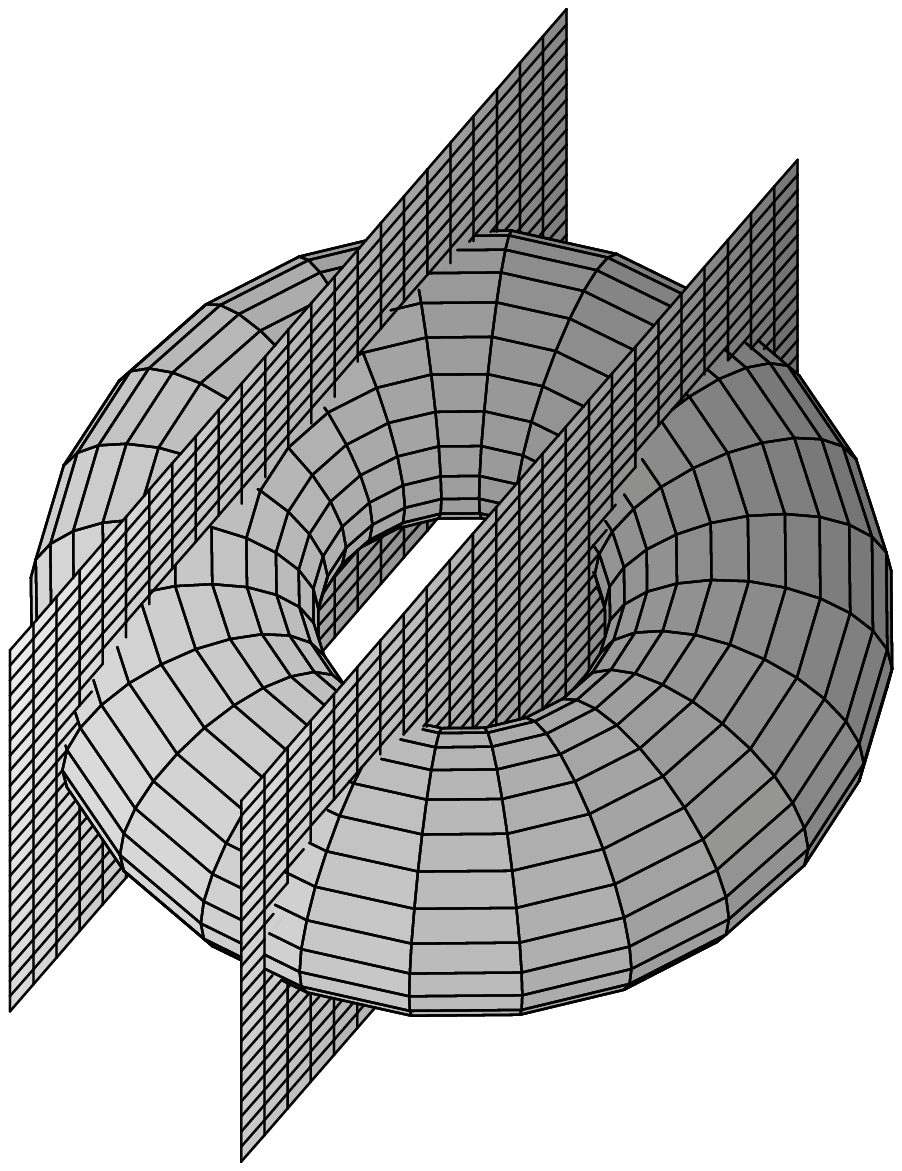}
  \epsfxsize=1.5in
  \epsffile{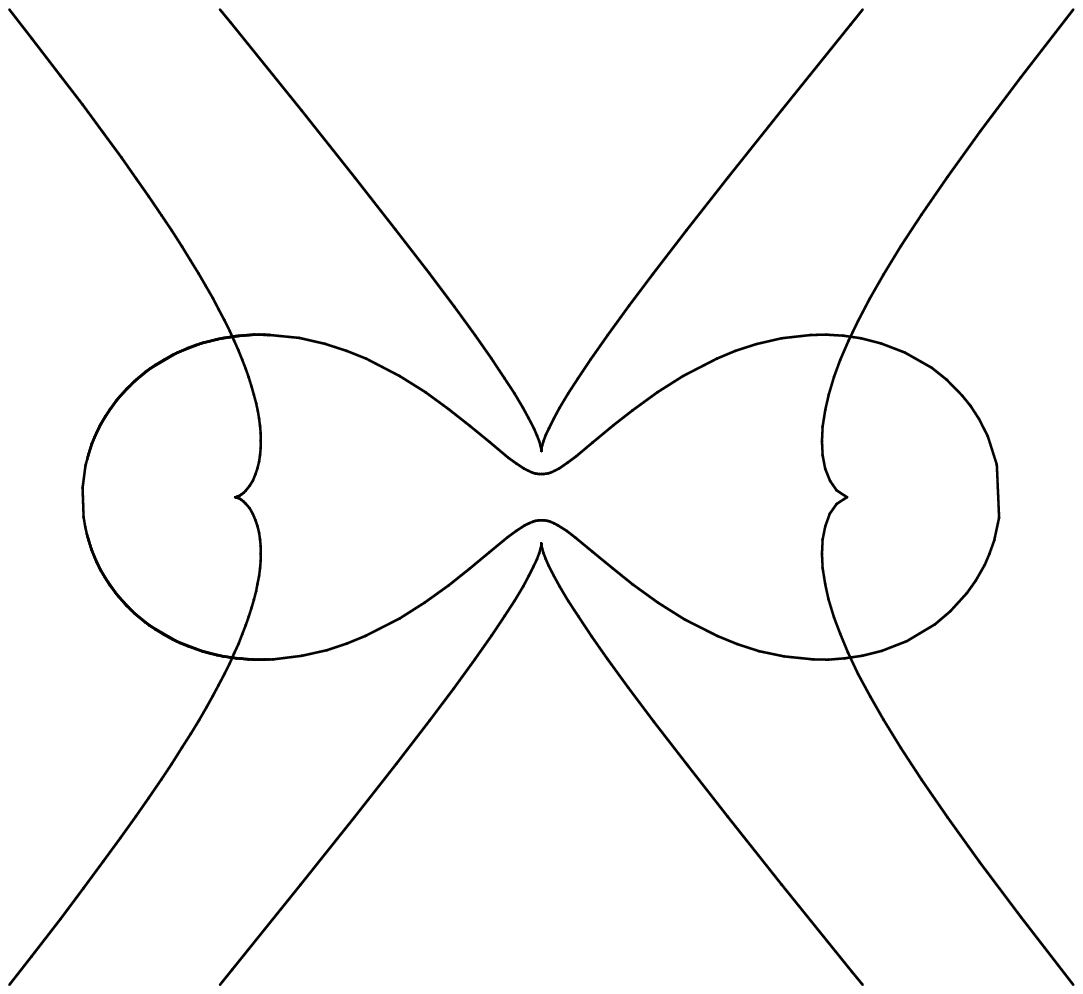}
\epsfxsize=1.5in
\epsffile{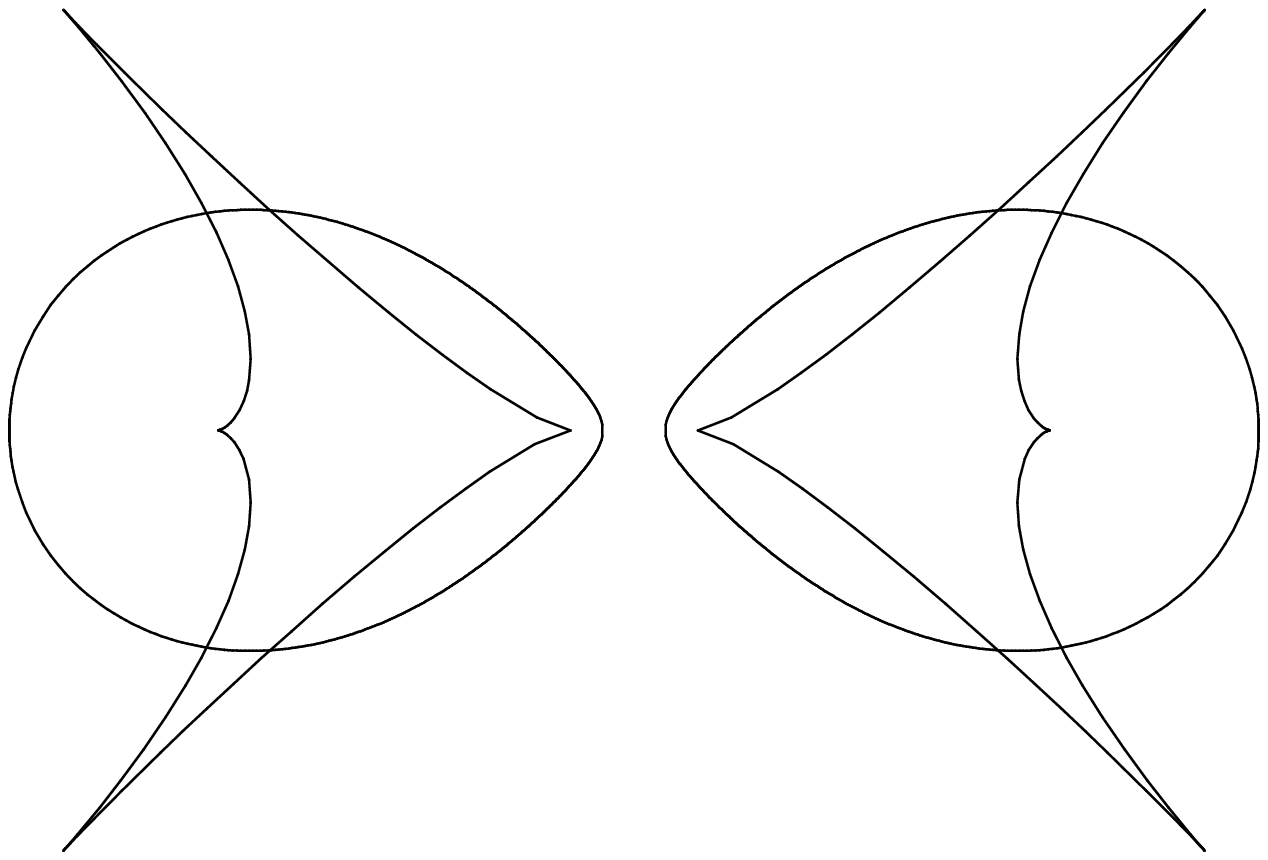}
  \end{center}
  \caption{\small Two plane sections of a torus close to a singular section, together with their evolutes. 
One connected component becomes two ovals, and four inflexions disappear while, locally,
two vertices become six.}
  \label{fig:torus}
  \end{figure}

Here is a simple example. Consider a round torus in 3-space, obtained by rotating a circle about an axis
in the plane of the circle but not intersecting it. This consists of elliptic and hyperbolic points,
separated by two circles of parabolic points along the `top' and `bottom' of the torus. (The
parabolic curves are far from generic but we shall stay clear of them.) We can take sections
by planes 
 parallel to the axis of rotation, as in Figure~\ref{fig:torus}. The sections  pass from a connected curve through a
nodal curve (a `figure eight') to two ovals. In the figure we have drawn the evolutes of the nonsingular sections: these have
cusps at the centres of curvature of the vertices. As the connected curve splits,
two vertices (maxima of curvature) come into coincidence at the crossing. 
 After the transition, when there are two components of the curve,
{\em three} vertices (one maximum and two minima of curvature) emerge from the crossing on each component.
 This transition, two local vertices becoming six local vertices, is written
`$1+1\leftrightarrow 3+3$'. As regards transitions on vertices this is one of two generic
situations at a hyperbolic point on a surface. However
as regards inflexions it is special since the figure-eight level curve itself has an inflexion on each branch at the crossing
point (in the terminology of \cite[p.282]{solid-shape} the crossing point is a flecnode for both asymptotic directions). This allows
a transition on the inflexions of the plane sections whereby $2+2\leftrightarrow 0+0$: two inflexions on each branch becomes
none. Taking into account both vertices and inflexions this becomes one of the types
`$\H_7$' below. If we take a hyperbolic point on the torus at which the tangent
plane is {\em not} parallel to the axis of rotation, then it can be shown that neither branch of the nodal curve has an inflexion. While
the transition on vertices remains as $1+1\leftrightarrow 3+3$, the inflexions become $2+0\leftrightarrow 1+1$ or
$1+1\lr 2+0$;
this is one of the `$\H_1$' cases in the notation below, which occur in regions of the surface.

\medskip

The paper is organized as follows.  In \S\ref{s:vertices-inflexions} we state the main results in the various
cases. In \S\ref{s:hyperbolic} we describe the various patterns which arise in the hyperbolic case, and the limiting
curvatures. The hyperbolic region of our surface $M$ is divided into subregions
according to the possible patterns, separated by a set which we call the
vertex transition (VT) set. This set consists of those hyperbolic points {\bf p} in $M$ for which one of 
the smooth local components of the intersection between $M$ and the tangent plane at {\bf p} has a vertex.
The VT set is difficult to calculate in particular cases
but in \S\ref{ss:VT} we give some  examples and
explain how the VT set approaches the parabolic curve on $M$.
In \S\ref{s:elliptic} we turn to the elliptic case, concentrating on umbilic points since the general
elliptic case is very simple. In \S\ref{s:parabolic} and \S\ref{s:cog} we cover the remaining
cases, parabolic point and cusp of Gauss respectively. Finally in \S\ref{s:conclusion} we summarize and
add some remarks on the material of the paper.

\smallskip\noindent\begin{small}
{\em Acknowledgements:} 
This work is a part of the DSSCV project supported by the IST Programme of the European Union (IST-2001-35443).
The first author was supported and the second author partially supported by this grant. We are also very grateful
to Terry Wall, Bill Bruce and Vladimir Zakalyukin for helpful suggestions.
\end{small}

\section{The vertex and inflexion sets}\label{s:vertices-inflexions}

We always assume that our surface is locally given by an equation $z=f(x,y)$ for some smooth function $f$, with
the tangent plane at the origin given by $z=0$. Thus our family of curves is $f(x,y)=k$ for constants
$k$ close to 0, and $(x,y)$ close to $(0,0)$. 
(In some cases the set $f(x,y)=k$ is non-empty only for one sign of $k$.)
We also take the $x$ and $y$ axes to be in principal directions at the origin, so that the surface $M$
assumes the local Monge form
\begin{eqnarray}
f(x,y)&=&\textstyle{\frac{1}{2}}\displaystyle (\kappa_1x^2 + \kappa_2y^2)+b_0x^3+b_1x^2y+b_2xy^2+b_3y^3\nonumber\\
&+&c_0x^4+
c_1x^3y+c_2x^2y^2+c_3xy^3+c_4y^4\nonumber \\
&+&d_0x^5+d_1x^4y+
d_2x^3y^2+d_3x^2y^3+d_4xy^4+d_5y^5+ \mbox{ h.o.t}.  
\label{eq:monge}
\end{eqnarray} where $\kappa_1,\kappa_2$ are the principal curvatures at $\pp$.  We often scale the surface
(multiply $x, y$ and $z$ by the same nonzero constant) so that $\kappa_1=2$ and the coefficient of
$x^2$ is therefore 1.

We use subscripts to denote partial derivatives: $f_x=\frac{\partial f}{\partial x}$ etc.

To such  $f$ we assign two functions $V_f$ and $I_f$  whose zero-level sets $V_f=0$  and  $I_f=0$ are respectively the 
sets of all vertices and inflexions of the
plane curves $f(x,y)=k$ for constants $k$.  We shall consider both the `vertex 
function' $V_f$ and the `vertex set' $V_f=0$. In fact for each of the generic cases of
elliptic, hyperbolic, parabolic and cusp of Gauss points of $M$ we shall go through the following steps.

\medskip\noindent
 $\bullet$ \ Calculate $V_f=0$ and $f=0$ and their Taylor expansions at the origin. In each 
case there will be several branches, some
of which may be singular.  (The same also applies to $I_f$.)
 \\
$\bullet$ \ Decide the possible relative positions of the branches of $f=0$ and $V_f=0$ (and $I_f=0$). These
can be indicated on diagrams.
\\
$\bullet$ For $k$ small, the level sets $f=k$  are close to the zero level set $f=0$. We can read off the pattern of 
vertices (and inflexions) from the diagrams above.
\\
$\bullet$  Calculate the limiting curvature at vertices of the section $f=k$, when $k\to 0$.

\medskip

To obtain the function $V_f$ we argue as follows. We want to find the vertices on a {\em smooth} curve $f(x,y)=k$. 
For this purpose we may assume that
locally the curve is given by $y=h(x)$ for a smooth $h$, that is $f(x,h(x))=k$ is an identity. Then the vertex
 condition is simply
$\kappa'(x)=0$ where  $\kappa (x) = \frac{h''(x)}{(1+(h'(x))^2)^{3/2}}$ 
is the curvature of $y=h(x)$. Working out the derivatives of $h$ in terms of those of $f$ and
clearing denominators we arrive at the following. The vertices of any smooth curve $f(x,y)=k$ will
be at the intersections with the set $V_f=0$, where
 \begin{eqnarray}
V_f &=&
(f_x^2+f_y^2)(-f_y^3f_{xxx}+3f_xf_y^2f_{xxy} -3f_x^2f_yf_{xyy}
+f_x^3f_{yyy}) \nonumber \\
&&+3f_xf_y(f_y^2f_{xx}^2+(f_x^2+f_y^2)f_{xx}f_{yy}-f_x^2f_{yy}^2) \nonumber \\
&&+6f_xf_yf_{xy}^2(f_x^2-f_y^2) \nonumber \\
&&+3f_{xy}(f_{xx}f_y^4-3f_x^2f_y^2(f_{xx}-f_{yy})-f_{yy}f_x^4) .
\label{eq:vertex}
\end{eqnarray}
The square of the curvature, $\kappa^2$, of the curve $f(x,y)=k$ at $(x,y)$ is
\begin{equation}
\kappa^2=\frac{(f_{xx}f_y^2-2f_{xy}f_xf_y+f_{yy}f_x^2)^2}{(f_x^2+f_y^2)^3},
\label{eq:kappa2}
\end{equation}
so that the {\em inflexion condition} is $I_f=0$ where
\begin{equation}
I_f(x,y)= f_{xx}f_y^2-2f_{xy}f_xf_y+f_{yy}f_x^2
\label{eq:infl}
\end{equation}
is the usual {\em Hessian determinant} of $f$.

The following result gives the number of intersections of the level set $f(x,y)=k$ with $V_f=0$ and $I_f=0$, as $k$ passes
through 0.  

\begin{theorem} \label{vertices-inflexions} Let $f=k$ be a section of a generic surface $M$ by a plane close to the
 tangent plane at  {\bf p}, $ k=0$ corresponding
with the tangent plane itself. Then for every sufficiently small 
open neighbourhood  $U$ of {\bf p} in $M$,  there exists $\varepsilon>0$ such that $f=k$ has exactly $v(\pp)$ vertices and $i(\pp)$ 
inflexions lying in $U$, for every $0<|k|\leq\varepsilon$, where $v(\pp)$ and $i(\pp)$ satisfy the following equalities. 
We also use $\leftrightarrow$ to indicate the numbers of vertices or inflexions on either
side of a transition, local to the singular point on $f=0$, when $f=k$ has two branches. The notation $m+n$ indicates the numbers 
of vertices or inflexions on the two branches.
\begin{enumerate}
\item[{\rm (E)}] If {\bf p} is an elliptic point, then for one sign of $k$ the section is locally empty; in the non-umbilic case,
for the sign of $k$ yielding a locally nonempty intersection we have $v(\pp)=4$, $i(\pp)=0$.
Likewise if {\bf p} is a generic\footnote{The genericity
assumption can be stated explicitly: the quadratic terms of $f$ should not divide the
cubic terms. See \S\ref{s:elliptic}.} umbilic point, then $v(\pp)=6$,  $i(\pp)=0$.
(This is already well-known: see for example~\cite[\S15.3]{porteous}.)

\item[{\rm (H)}] If {\bf p} is a hyperbolic point v(\pp) satisfies
one of the following. 

For {\bf p} lying in open regions of $M$ we have \\   $2+2 \leftrightarrow 2+2$ or  $1+1 \leftrightarrow 3+3$. \\
In other cases, occurring along curves or at isolated points of $M$, we can have in addition\\
$3+2\lr 2+1$ or $3+1\lr 2+2$.\\
See \S\ref{ss:hyp-vertex-inflexion} for an explanation of the different cases.

\smallskip\noindent
Also using the same notation, $i(\pp)$ satisfies:  $1+1 \leftrightarrow 0+2$ or  $1+2 \leftrightarrow 0+1$; the
full list is in Table~\ref{table:hyper2}.

\item[{\rm (P)}]
If {\bf p} is a parabolic point but not a cusp of Gauss, $v(\pp)=3$,  $i(\pp)=2$.

\item[{\rm (ECG)}] If {\bf p} is an elliptic cusp of Gauss, $v(\pp)=4$ ,  $i(\pp)=2$ for one sign of  $k$, and the
level set is empty for the other.
\item[{\rm (HCG)}]
If {\bf p} is a hyperbolic cusp of Gauss, we have:\\
$v(\pp) : 1+3\leftrightarrow 4+4$ or $ 2+2\leftrightarrow 4+4, $ and for each of these, we can have any of\\
$ i(\pp): 1+1\lr 0+0$ or $2+2\leftrightarrow 0+2$ or $1+1\lr 0+4$
\end{enumerate}
\end{theorem}

We split the proof of Theorem~\ref{vertices-inflexions} into different cases discussed in 
the relevant sections, in which we also carry out a closer investigation
of the geometry of the sets $V_f=0$ and $I_f=0$.

\section{Hyperbolic case}\label{s:hyperbolic}
Recall that at a hyperbolic point {\bf p} of  a surface, the principal curvatures $\kappa_1,\kappa_2$ are 
not zero and have opposite signs. After scaling, 
$f$ can be taken in (\ref{eq:monge}) to have quadratic part $x^2-a^2y^2$ where $a>0$. 
We shall write $V_h$ for $V_f$ in this case, and likewise $I_h$ for $I_f$.

\subsection{Patterns of vertices and inflexions on the level sets}\label{ss:hyp-vertex-inflexion}

\begin{proposition}\label{prop:branches-hyperbolic}
{\rm (i)} \ The vertex set $V_{h}=0$ has exactly four smooth branches $VH_1, VH_{2},VH_{3}$, $VH_{4}$ through $(0,0)$, where
$VH_1$ is tangent to the principal direction   $ x=0 $,
$VH_{2}$ is tangent to the principal direction $ y=0 $,
$VH_{3}$ is tangent to the asymptotic direction $x - ay =0 $ and
$VH_{4}$ is tangent to the asymptotic direction  $ x + ay =0 $.

\smallskip\noindent
{\rm (ii)} \ The level sets $f=0$  and $I_{h}=0$   have exactly two smooth branches in a neighbourhood of $(0,0)$, 
one of them being tangent to $ x-ay=0$ and the other one to $ x + ay =0 $.
\label{prop:vert-inf}
\end{proposition}
The proof for $V_h$ can be done in several ways. We can use the technique exemplified
in \S\ref{ss:proof-para-vert-inflex}, that is, blowing up combined with the implicit function theorem,
or, in the present case, we can even prove that $V_h$ is ${\mathcal R}$-equivalent as a function
to its lowest terms, which are $192a^4(1+a^2)xy(x - ay)(x+ay)$. The functions $f$ and $I_h$ are
Morse functions, hence equivalent to their quadratic parts.

In order to verify the conclusions of Theorem~\ref{vertices-inflexions} in the hyperbolic case we need to determine the relative 
positions of the branches of $f=0$ and $V_f=0$ (and $I_f=0$)
which are tangent to one another at the origin.  
To do this we need the higher terms of the Taylor expansions of those branches
with the same tangents. The branches $VH_1$ and $VH_2$
present no problems since they are always transverse to the branches of the
level set $f=0$. For the branches $VH_3$ and $VH_4$ we use Proposition~\ref{prop:branches-hyperbolic} and
substitute for example $x=ay+x_2y^2+x_3y^3+\mbox{higher terms}$ into the expression the vertex set 
$V_{h}$, for the branch $VH_{3}$.

\medskip\noindent
{\bf Notation} \ Certain expressions occur often in our formulae so we introduce some notation for them.\\
$f^{(n)}(a)$ means the result of substituting $x=a, \ y=1$ in the homogeneous part of degree $n$ in the Taylor expansion
of $f$. (We write this rather than the more precise $f^{(n)}(a,1)$.) For example,\\
$f^{(3)}(a)=b_0a^3+b_1a^2+b_2a+b_3$, and \\
$f^{(4)}(a)=c_0a^4+c_1a^3+c_2a^2+c_3a+c_4$.\\
Similarly $f^{(n)}(-a)$ is the result of
substituting $x=-a, y=1$ in the same homogeneous polynomial of degree $n$.

\begin{proposition}\label{pro:hyp-vertices}
{\rm (i)} \ The branches $VH_{3}, \ VH_{4}$ of the vertex set have the following 3-jets:
\begin{eqnarray*}
 &&VH_{3}: \\
x &=& ay  -\frac{1}{2a}f^{(3)}(a)y^2 \\
&+& \left.\frac{1}{4a^3(1+a^2)}\right( f^{(3)}(a)(3b_0a^5+b_1a^4+(5b_0-b_2)a^3
+(3b_1-3b_3)a^2+b_2a-b_3)\\
&&\hspace*{1in} - \ 4a^2(1+a^2)f^{(4)}(a) \left)\phantom{\frac{1}{2}}\right. \!\!\!\!\!\! y^3 \\
&=& ay +x^+_{2v}y^2+x^+_{3v}y^3 \ \mbox{say},\\
VH_{4}: x&=& -ay +x^-_{2v}y^2+x^-_{3v}y^3
 \mbox{(obtained by replacing $a$ with $-a$ in the above.)}
\end{eqnarray*}

\medskip\noindent
{\rm (ii)} \ The branches of $f=0$ have the following 3-jets:
\begin{eqnarray*}
x &=& ay   -\frac{1}{2a}f^{(3)}(a)y^2\\
&+&
\left.\frac{1}{8a^3}\right( f^{(3)}(a)(5b_0a^3+3b_1a^2+b_2a-b_3) -4a^2f^{(4)}(a)\left)\phantom{\frac{1}{2}}\right. \!\!\!\!\!\! y^3\\
&=& ay+x^+_{2v}y^2+x^+_{3f}y^3 \ \mbox{say, and}\\
x&=& -ay+x^-_{2v}y^2+x^-_{3f}y^3 \mbox{
 obtained by replacing $a$ with $-a$ in the above}.
\end{eqnarray*}
\label{prop:vert-taylor}
\end{proposition}

\vspace*{-0.5cm}

It is evident from (i) and (ii) of the above Proposition that the branches of vertex set and those
of the curve $f=0$ have at least 3-point contact at the origin: their Taylor expansions
agree up to order two. This also means that they have the same osculating circle (circle of curvature)
at the origin. The condition for them to have (at least) 4-point contact is that the
terms in $y^3$ agree also.
 After some manipulation, this 4-point contact condition comes to
the following.

\begin{proposition}\label{prop:4-pointcontactcondition}
{\bf Four-point contact condition} \ 
The condition for the vertex branch $VH_3$ to have (at least) 4-point contact with
the corresponding branch of $f=0$ at the origin is
\begin{equation}
f^{(3)}(a)(b_0a^5-b_1a^4+(5b_0-3b_2)a^3-(5b_3-3b_1)a^2+b_2a-b_3)
-4a^2(1+a^2)f^{(4)}(a)=0.
\label{eq:VT1}
\end{equation}
The condition for $VH_4$ to have (at least) 4-point contact with
the corresponding branch of $f=0$ is obtained by replacing $a$ by $-a$:
\begin{equation}
-f^{(3)}(-a)(b_0a^5+b_1a^4+(5b_0-3b_2)a^3+(5b_3-3b_1)a^2+b_2a+b_3)
-4a^2(1+a^2)f^{(4)}(-a)=0.
\label{eq:VT2}
\end{equation}
\label{prop:VT} 
\end{proposition}

\vspace*{-1cm}

For a generic surface $M$, (\ref{eq:VT1}) or ({\ref{eq:VT2}) then imposes one condition
on the point {\bf p} and can therefore be expected to hold for points
{\bf p} along  one or more {\em curves} on $M$. We call this the {\em vertex transition set} (VT set) on $M$.

\begin{remarks}

\noindent
{\rm
{\bf (1)} \ The apparently rather complicated conditions in the above proposition actually state that {\em one or other of the branches of
the curve $f=0$ itself---the intersection between the surface $M: z=f(x,y)$ and its tangent plane---has a vertex}. In fact we
have the general result:

\begin{quote} {\em For any $f$ giving a hyperbolic point at the origin,
a branch of the curve $f=0$ and the corresponding branch of the vertex set have the same
order of contact with their common osculating circle.}
\end{quote}

Thus at a point of the VT set, the corresponding branches of $f=0$ and of the vertex set
both have vertices. We shall not use this fact here, but discuss the result and its
consequences elsewhere.

\medskip\noindent

{\bf (2)} \ Note in particular that (\ref{eq:VT1}) holds if $x-ay$ is a factor
of both the cubic and quartic terms of the expansion of $f$. This is a {\em biflecnode}
in the terminology of Koenderink \cite[p.296]{solid-shape}. As a special case, one of
(\ref{eq:VT1}), (\ref{eq:VT2}) will hold at every point of a ruled surface, since the whole
line in one of the asymptotic directions lies on the surface. From the point of view
of the VT set, both ruled surfaces and surfaces of revolution (see \S\ref{ss:VT}) are highly non-generic.
}
\end{remarks}

Note that the 4-point contact condition can be regarded as a formula for $f^{(4)}(\pm a)$ in terms of the
lower degree coefficients of the expansion of $f$. It can therefore be regarded as a formula for any of the
degree 4 coefficients $c_i$ in terms of the other $c_j$ and lower degree coefficients of $f$. In a similar
way we can write down the additional condition for $VH_3$ or $VH_4$ and the corresponding branch of $f=0$ to have
5-point contact. This can be written in the form $f^{(5)}(\pm a)=$ a polynomial in the lower degree coefficients,
but it is complicated and we shall not display it here. (As noted above, this is equivalent to the branch
of $f=0$, or of the vertex set, having a higher vertex.)

Analysing in a similar way the Taylor expansions of the inflexion function $I_h$  we find the
following.
\begin{proposition}\label{prop:hyp-inflexions}
The branches of the inflexion curve $I_h=0$ have the following 3-jets:
\begin{eqnarray*}
x &=& ay +\left.\frac{1}{8a^3}
\right(-3f^{(3)}(a)(3b_0a^3+b_1a^2-b_2a-3b_3)+8a^2f^{(4)}(a)\left)\phantom{\frac{1}{2}}\right. \!\!\!\!\!\! y^3 \\
&=& ay+x^+_{3i}y^3 \ \mbox{say},
\\
 \mbox{and} &&\\
x&=& -ay+x^-_{3i}y^3 \  \mbox{obtained by replacing $a$ with $-a$ in the above}.
\end{eqnarray*}
\label{prop:inf-taylor}
\end{proposition}

\vspace*{-0.5cm}

Note that there are no quadratic terms in these expansions: the branches of the inflexion curve $I_{h}=0$
{\em themselves have inflexions} at the origin. Accordingly the branches of $I_h=0$ and
$f=0$ tangent to $x=ay$, say, have 2-point contact unless the branch of $f=0$
also has an inflexion (that is, $f^{(3)}(a)=0$).

Altogether the possibilities for contact between branches of the vertex and inflexion sets
and the branches of $f=0$ in the present hyperbolic cases are as follows.

\medskip\noindent
{\bf Notation} \\
$\V_1$ \ A branch of $f=0$ and of $V_h=0$ have the minimum 3-point contact,\\
$\V_2$ \ A branch of $f=0$ and of $V_h=0$  have 4-point contact; see (\ref{eq:VT1}) or (\ref{eq:VT2}), \\
$\V_3$ \ A branch of $f=0$ and of $V_h=0$ have 5-point contact,

\smallskip\noindent
$\I_1$ \ A branch of $f=0$ and of $I_h=0$ have the minimum 2-point contact, \\
$\I_2$ \ A branch of $f=0$ and of $I_h=0$ have 3-point contact. \\ 
$\I_3$ \ A branch of $f=0$ and of $I_h=0$ have 4-point contact.

The possible ways of combining these at the two branches of $f=0$ tangent to $x=\pm ay$ are therefore as shown
in Table~\ref{table:hyper1}. Here
`codim' refers to the codimension of the locus of these points in the hyperbolic region.

\medskip

\begin{table}
\begin{center}
\begin{tabular}{|c|c|c|c|l|}\hline
Symbol&$x=ay$ branch & $x=-ay$ branch & `codim' & Comment  \\  \hline\hline
$\H_1$&$\V_1\I_1$ & $\V_1\I_1$ & 0 & the most generic case \\ \hline
$\H_2$&$\V_2\I_1$ & $\V_1\I_1$ & 1 &  along curves in the VT set \\ \hline
$\H_3$&$\V_2\I_1$ & $\V_2\I_1$ & 2 & self-intersections of the VT set \\ \hline
$\H_4$&$\V_3\I_1$ & $\V_1\I_1$ &  2 & isolated points of the VT set \\ \hline
$\H_5$&$\V_1\I_2$ &  $\V_1\I_1$ & 1 & curves in the hyperbolic region\\ \hline
$\H_6$&$\V_2\I_1$ &  $\V_1\I_2$ & 2 &  isolated points \\ \hline
$\H_7$&$\V_1\I_2$ & $ \V_1\I_2$ & 2 &  isolated points \\ \hline
$\H_8$&$\V_2\I_3$ &  $\V_1\I_1$ & 2 & $\V_2\I_2$ and $\V_1\I_3$ do not occur \\ \hline
\end{tabular}
\end{center}
\caption{\small The possibilities for contact between $f=0$ and the vertex and inflexion curves. See Propositions~\ref{prop:VT}
and \ref{prop:H1}, and Lemma~\ref{lemma:VI}
for further information.}
\label{table:hyper1}
\end{table}

\medskip

For the `most generic' case $\H_1$, we give in Figure~\ref{fig:hyp-geom-vertex-inflex} the three possible ways
(up to rotation or reflexion of the diagram) in which the different elements can intersect.  We use the
notation $2+2\leftrightarrow 2+2$ and $1+1\leftrightarrow 3+3$ to indicate the numbers
of vertices on the pair of branches of $f=k$ for small $k$ first of one sign and then of the other.
The inflexions in the first case follow the pattern $2+0 \lr 1+1$ whereas in the second case
the two patterns $2+0\lr 1+1$ and $1+1\lr 2+0$ occur. Examining  cases we find the following.

\begin{proposition}
{\rm (i)}  \ In the case $\H_1$, the vertex transition $2+2\leftrightarrow 2+2$ occurs when the left hand sides of
(\ref{eq:VT1}) and (\ref{eq:VT2}) have opposite signs. \\
{\rm (ii)}
The vertex transition $1+1\leftrightarrow 3+3$ occurs when the left hand sides of
(\ref{eq:VT1}) and (\ref{eq:VT2}) have the same sign.\\
{\rm (iia)} \ The inflexion transition $1+1\lr 2+0$ occurs when, in addition to {\rm (ii)}, `left hand sides both $<0$' is accompanied
by $f^{(3)}(a)f^{(3)}(-a)>0$ and `left hand sides both $>0$' by $f^{(3)}(a)f^{(3)}(-a)<0$.\\
{\rm (iib)} \ The inflexion transition $2+0\lr 1+1$ occurs when, in addition to {\rm (ii)}, `left hand sides both $<0$' is accompanied
by $f^{(3)}(a)f^{(3)}(-a)<0$ and `left hand sides both $>0$' by $f^{(3)}(a)f^{(3)}(-a)>0$.\\
\label{prop:H1}
\end{proposition}

The other cases also require an analysis of the order of branches of $f=0, V_h=0, I_h=0$
around each of the lines $x=\pm ay$. In Figure~\ref{fig:I2andV2} the principal cases for a single
branch tangent to $x=ay$ are illustrated. By putting these together with similar information at $x=-ay$,
and including the other branches of $V_h=0$ tangent to the two coordinate axes (see Proposition~\ref{prop:vert-inf})
we arrive at the classification in Table~\ref{table:hyper2}.

Here is an indication of the calculations which allow us to draw the cases in Figure~\ref{fig:I2andV2}.

\begin{lemma}
{\rm (i)} \  The condition for $\V_1\I_2$ on the branch tangent to $x=ay$ is $f^{(3)}(a)=0, f^{(4)}(a)\ne 0$
 and the configuration of $f=0, V_h=0, and I_h=0$
is determined by the sign of $f^{(4)}(a)$, as in Figure~\ref{fig:I2andV2}. \\
{\rm (ii)} \ The conditions for $\V_2\I_2$ or $\V_1\I_3$ on the branch tangent to $x=ay$ are $f^{(3)}(a)=f^{(4)}(a)=0,
f^{(5)}(a)\ne 0$, and
this situation is  in fact $\V_2\I_3$. 
The configuration of $f=0, V_h=0$ and $I_h=0$ is determined by the sign of $f^{(5)}(a)$, as in Figure~\ref{fig:I2andV2}. 
\label{lemma:VI}
\end{lemma}
{\bf  Proof} \ For (i), note that we require the branches of $I_h=0$ and $f=0$ tangent to $x=ay$ to have
the same 2-jet, and using the formulae of Propositions~\ref{prop:vert-taylor} and \ref{prop:inf-taylor} this
requires $f^{(3)}(a)=0$. They have different 3-jets provided $f^{(4)}(a)\ne 0$, since the coefficients
of $y^3$ in the two Taylor series are then $\frac{1}{a}f^{(4)}(a)$ and $-\frac{1}{2a}f^{(4)}(a)$ respectively.
Since the coefficient of $y^3$ in the Taylor series of $V_h$ is $-\frac{1}{a}$ we find the two orderings
of the branches depicted in Figure~\ref{fig:I2andV2}.\\
For (ii) we use the same Propositions, noting that $\I_2$ together with $\V_2$ imply $f^{(3)}(a)=f^{(4)}=0$
which in turn imply that the 3-jets of the branches of $V_h=0$ and $f=0$ agree. Further calculations then show that
the coefficients of $y^4$ in the three branch expansions are $f=0: -\frac{1}{2a}f^{(5)}(a), \ V_h=0: -\frac{5}{2a}f^{(5)}(a)$
and $I_h=0: \frac{5}{2a}f^{(5)}(a)$ from which the results now follow.

\begin{table}
\begin{center}
\begin{tabular}{|rl|c|c|l|}\hline
\multicolumn{2}{|c|}{Symbol}& Vertex transitions & Inflexion transitions &  Comment  \\  \hline\hline
$\H_1$&(i)& $2+2\leftrightarrow 2+2$&$2+0 \lr 1+1$ &Figure~\ref{fig:hyp-geom-vertex-inflex}(i)\\ \hline
&(iia)&        $1+1\leftrightarrow 3+3$ &      $1+1\lr 2+0$   & Figure~\ref{fig:hyp-geom-vertex-inflex}(iia)\\ 
&(iib) && $2+0\lr 1+1$ &Figure~\ref{fig:hyp-geom-vertex-inflex}(iib) \\ \hline
$\H_2$&& $3+2\leftrightarrow 2+1$ & $1+1\lr 2+0$ &  \\ 
& & & $2+0\lr 1+1$  & \\
& & & $0+2\lr 1+1$  & \\
& & & $1+1\lr 0+2$ & \\ \hline
$\H_3$& & $3+1\lr 2+2$& $2+0\lr 1+1$& \\ 
&&&$1+1\lr 2+0$&\\
&&&$0+2\lr 1+1$ &\\ \hline
$\H_4$& & && As for $\H_1$ \\ \hline
$\H_5$&(i) & $2+2\lr 2+2$ & $1+2\lr 1+0$ & \\ \hline
&(ii)& $1+1\lr 3+3$ & $1+2 \lr 1+0$ & \\ \hline
$\H_6$& & $3+2\lr 2+1$ & $1+0\lr 2+1$ & \\ \hline
$\H_7$&(i) & $2+2 \lr 2+2$ & $1+1 \lr 1+1$ & \\ \hline
 &(ii) &  $1+1 \lr 3+3$ & $2+2 \lr 0+0$ & Torus, Figure~\ref{fig:torus} \\ \hline
$\H_8$ && $3+2\lr 2+1$ & $1+1\lr 0+2$ & \\ 
& & &$0+2\lr 1+1$ & \\ \hline
\end{tabular}
\end{center}
\caption{\small The transitions on vertices and inflexions in the hyperbolic case.}
\label{table:hyper2}
\end{table}

\medskip

\begin{figure}[h!]
  \begin{center}
  \leavevmode
\epsfysize=3in
\epsffile{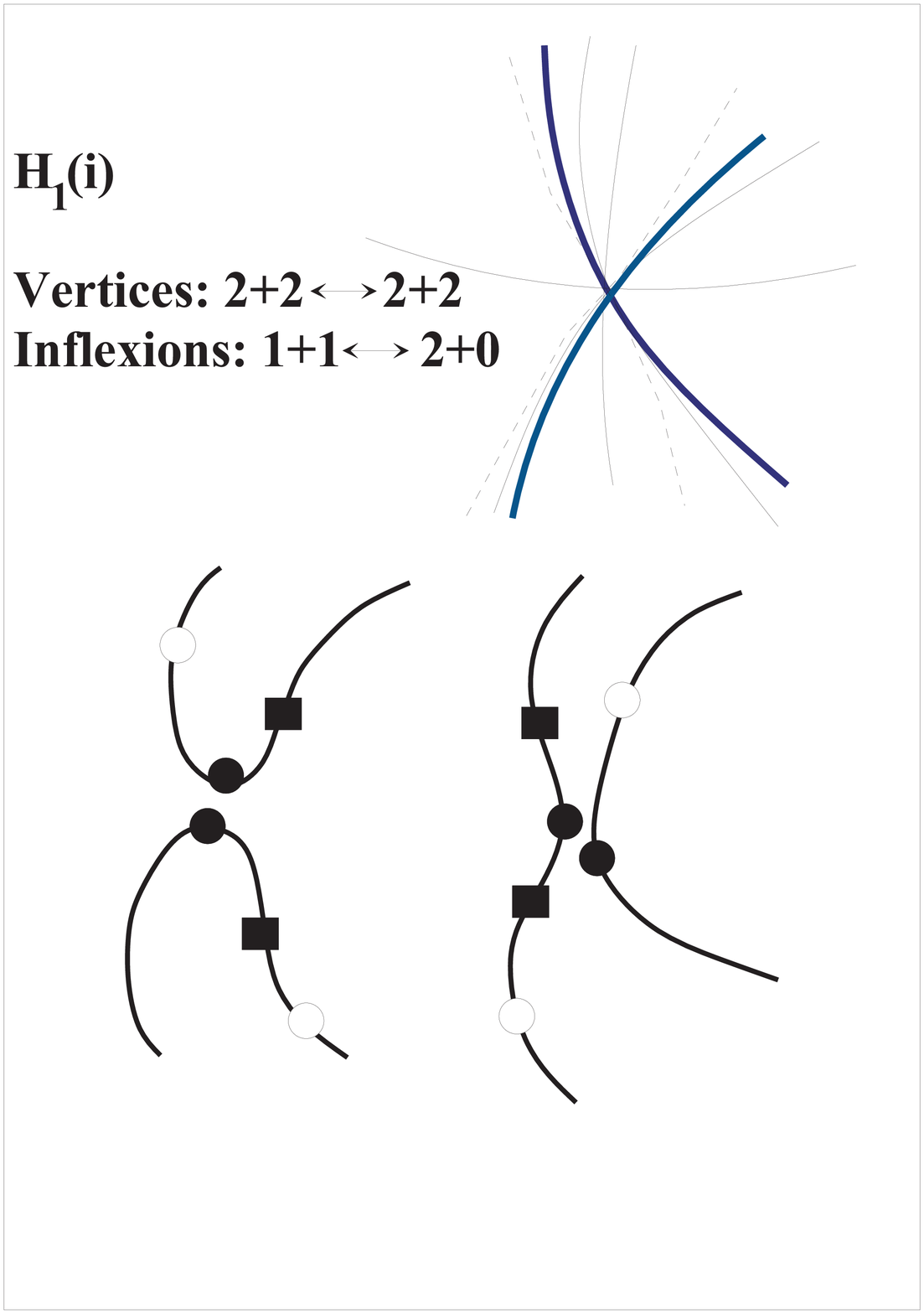}
\hspace*{0.35in}
\epsfysize=1.2in
\epsffile{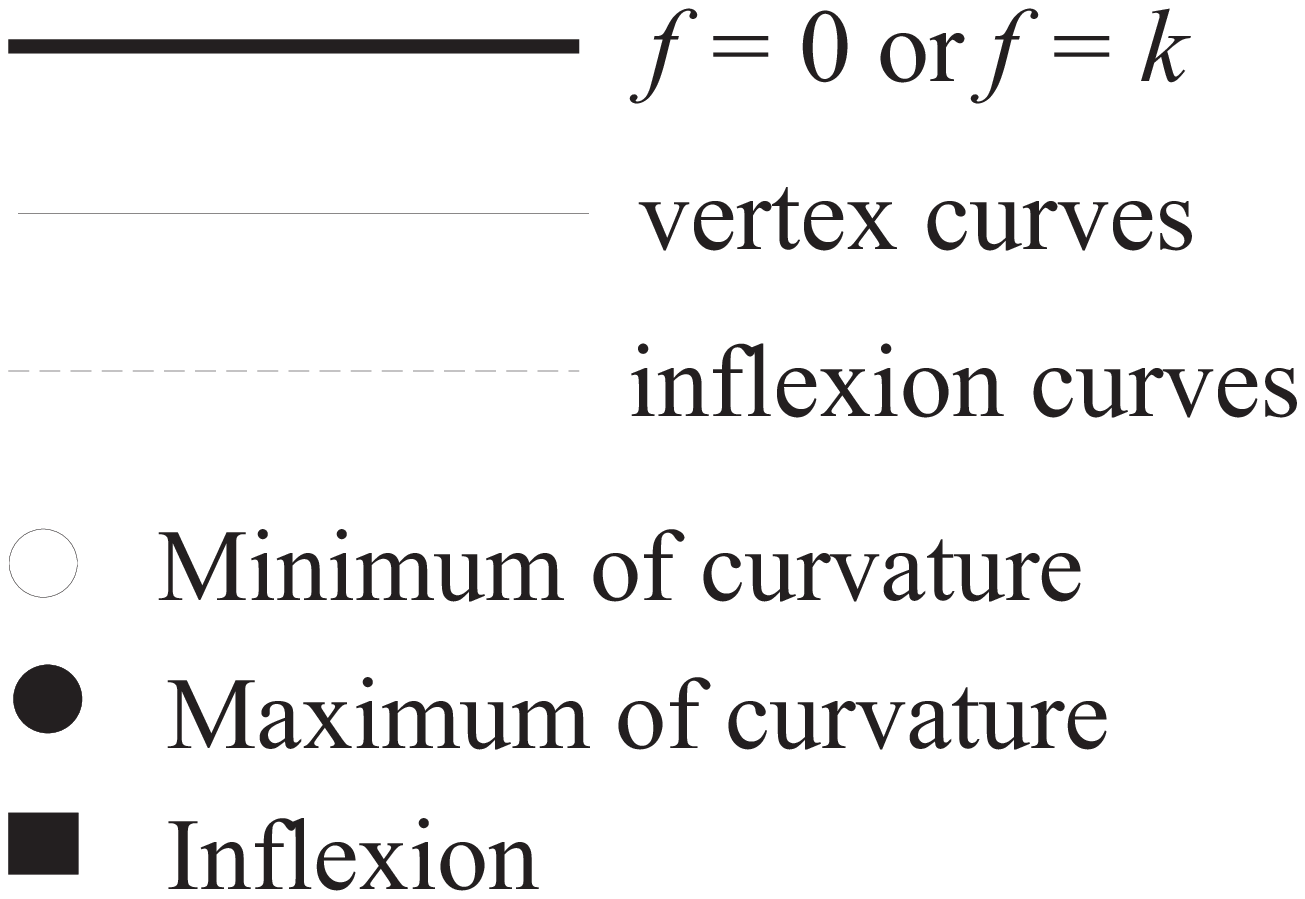}
\end{center}
\begin{center}
\epsfysize=3in
\epsffile{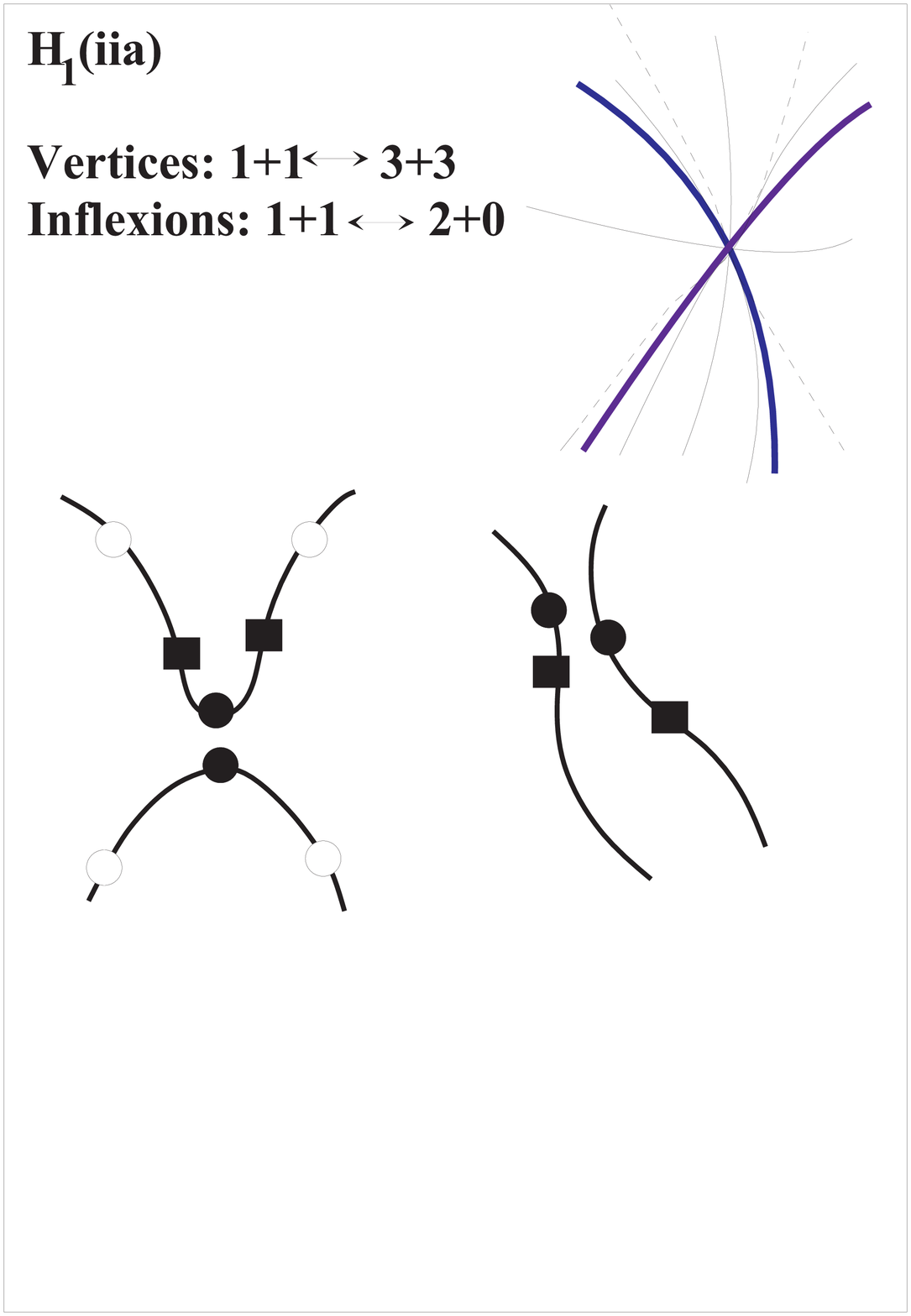}
\epsfysize=3in
\epsffile{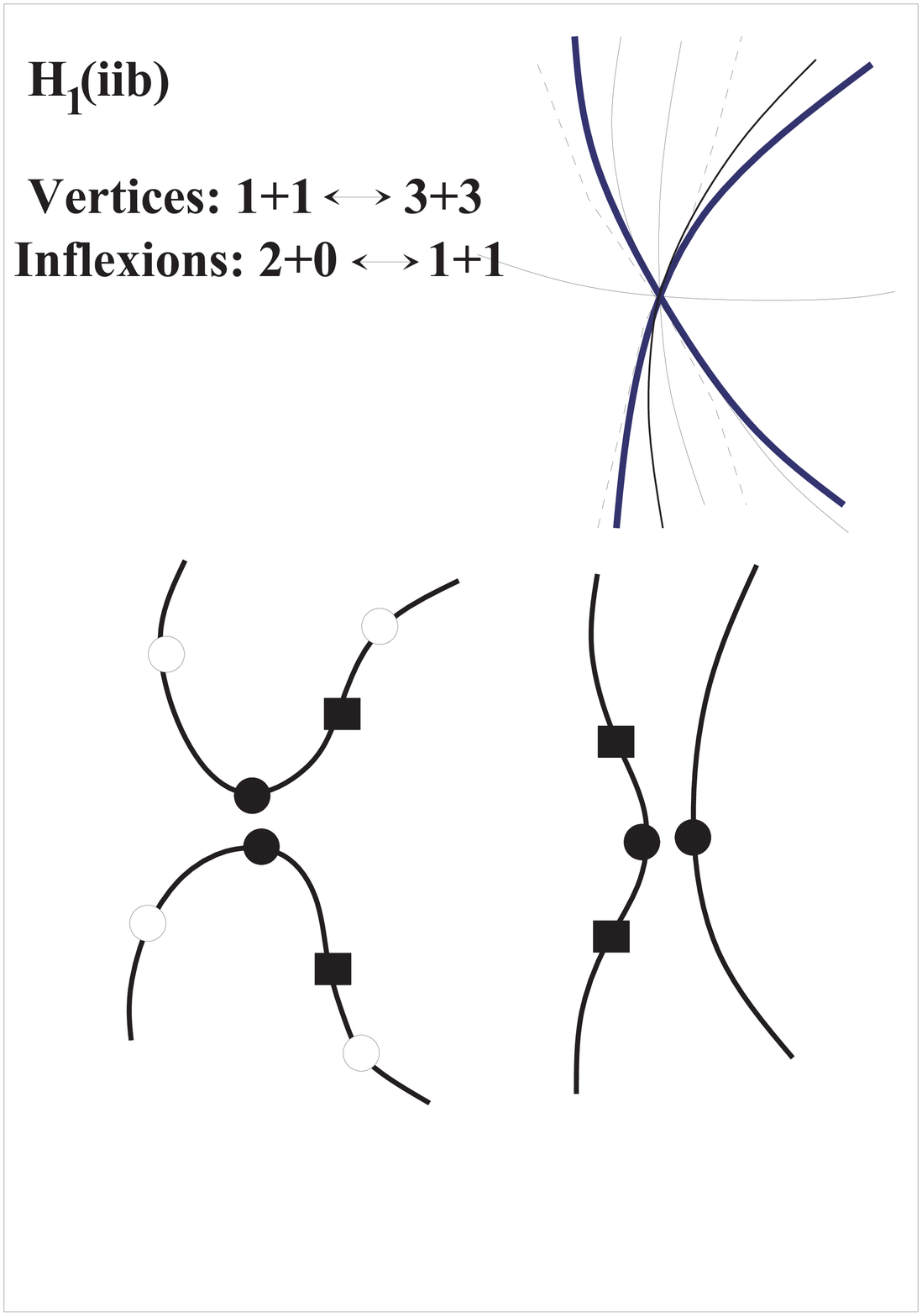}
\end{center}
  \caption{\small Arrangements of vertices and inflexions on the level sets of $f$, hyperbolic case ${\bf H}_1$ (see Table~\ref{table:hyper2}).
In each case, we show, above, the vertex and inflexion curves---that is, the loci of vertices and inflexions on the level sets of $f$---and,
below, a sketch of the level curves for $f<0$, $f>0$, showing the positions of these vertices and inflexions. The orientation chosen for
the branches of $f=k$ is shown in Figure~\ref{fig:growthkappa}.}
  \label{fig:hyp-geom-vertex-inflex}
  \end{figure}

\begin{figure}
  \begin{center}
  \leavevmode
\epsfysize=1.4in
\epsffile{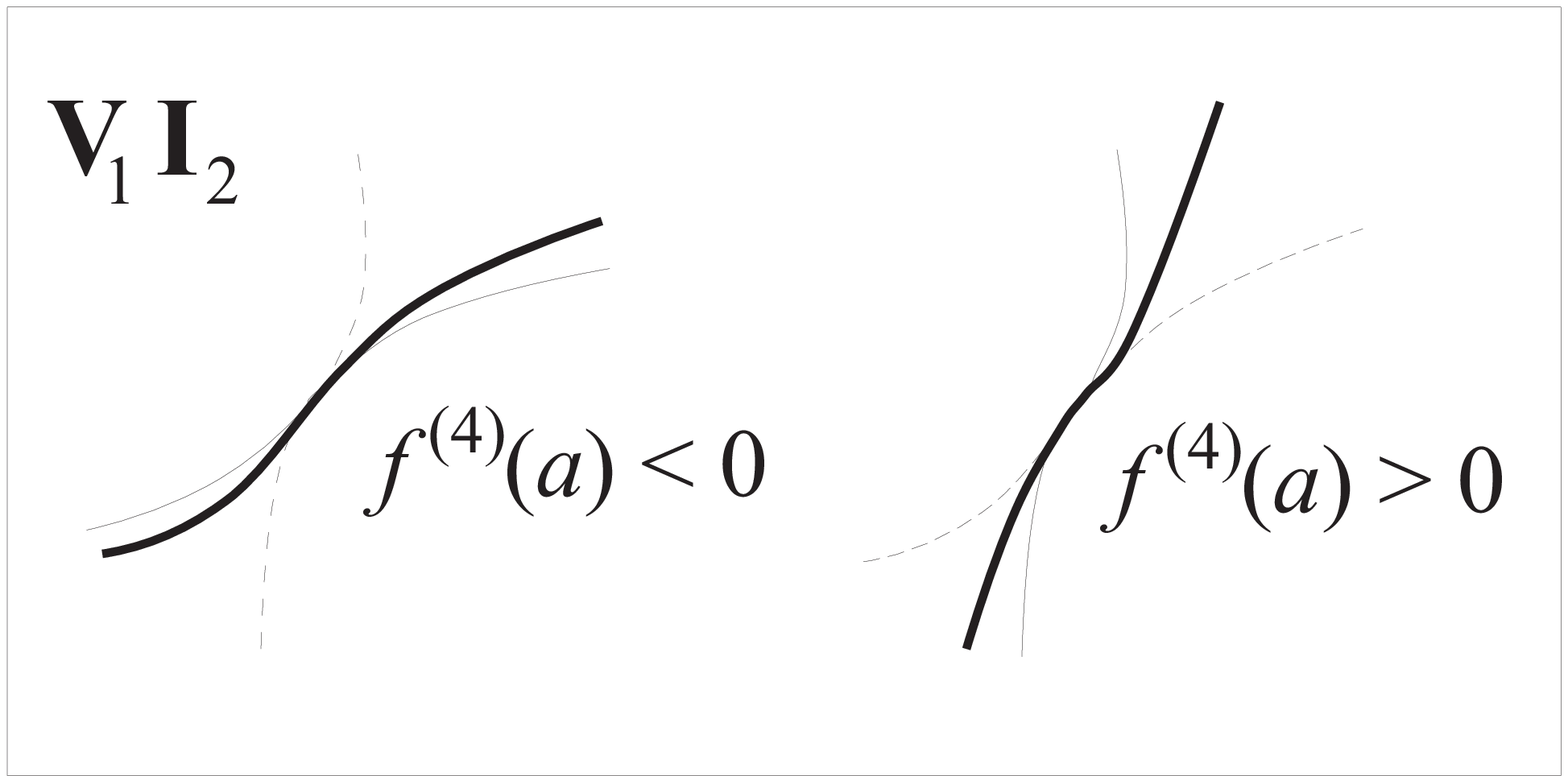}
\hspace*{0.35in}
\epsfysize=1.4in
\epsffile{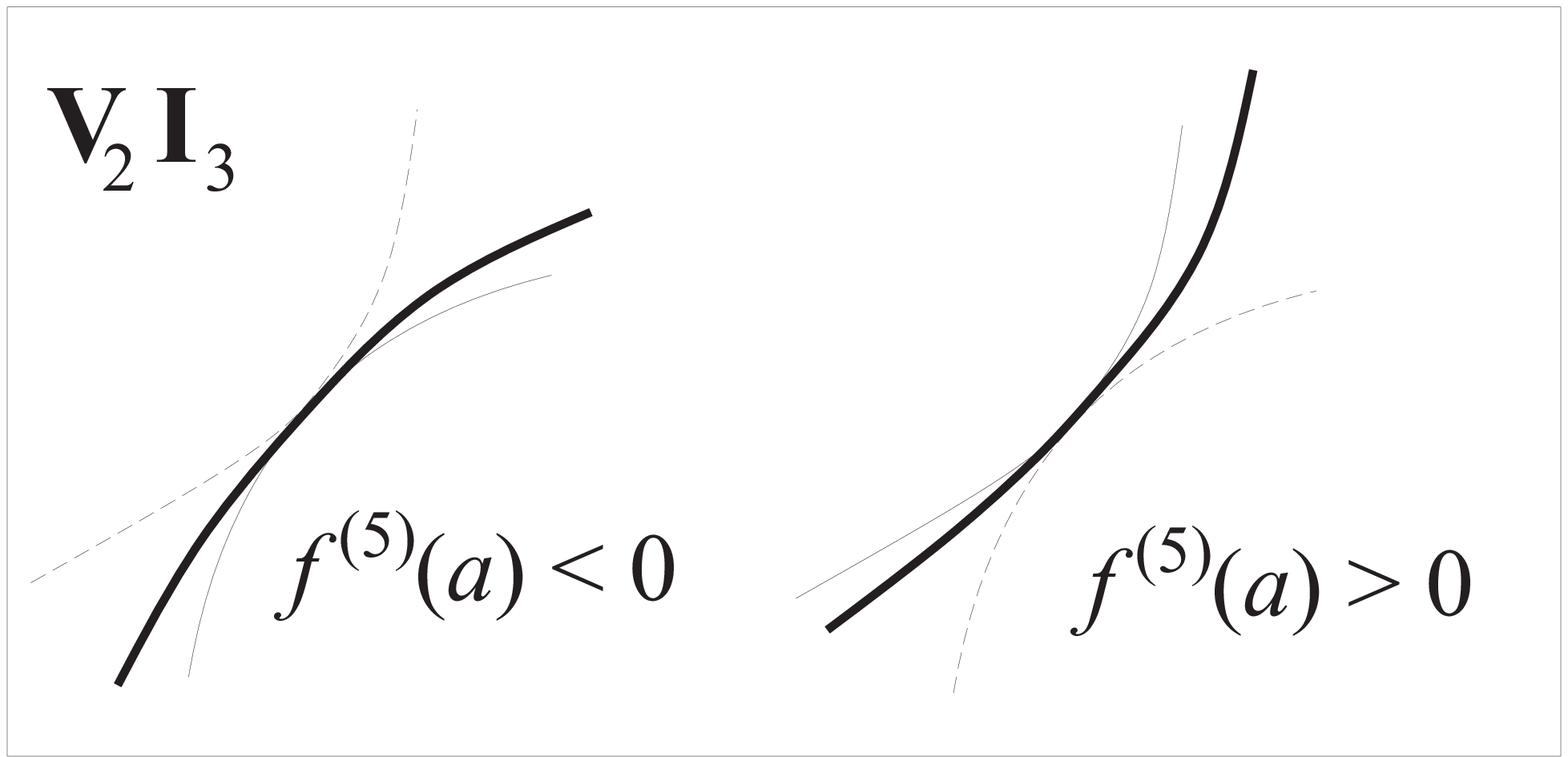}
\end{center}
\begin{center}
\epsfysize=1.15in
\epsffile{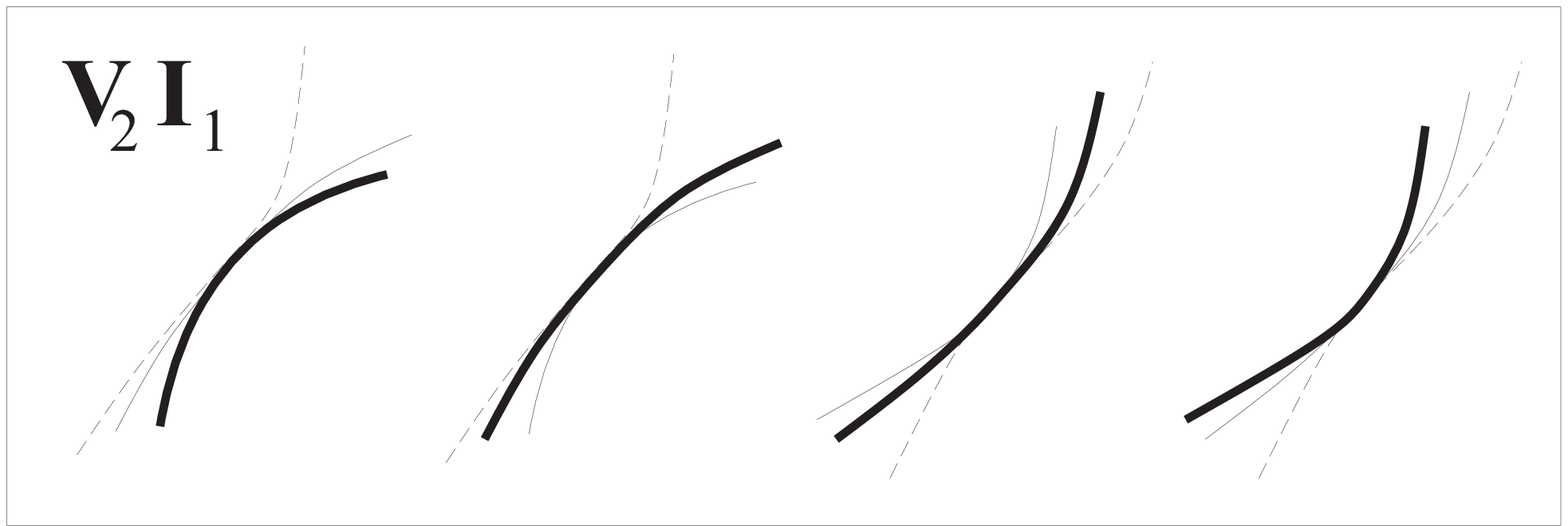}
\end{center}
  \caption{\small The arrangements of branches tangent to $x=ay$: thick line $f=0$, thin solid line the vertex curve $V_h=0$ and
dashed line the inflexion curve $I_h=0$. Three cases are illustrated here, the notation being that of Table~\ref{table:hyper1}.}
\label{fig:I2andV2}
  \end{figure}

\subsection{\bf Extrema of  curvature  and limiting curvature}

In order to  analyse the vertices further, we need to decide which vertices 
correspond to maxima and which to minima of curvature on the curve. (This is of significance when
we apply the results to the symmetry set and the medial axis, since only minima---indeed absolute minima---can contribute
to the latter.) We proceed as follows. The different 
branches of the vertex set locally
divide the plane into regions where the derivative  $\kappa'$ of the curvature $\kappa$
(with respect to any regular parametrisation 
of the curve) has a constant sign, the vertex branches 
being the loci of points where this derivative vanishes. Note that the sign of $\kappa'$ does {\em not} depend on the orientation of the curve.
However $\kappa'$ has the same sign as the vertex condition $V_h(x,y)$.

To decide the sign of $\kappa'$, for instance in the (local) region between the vertex branches tangent to  $y=0$ and $x-ay=0$, 
let us then check the sign of $V_h(x,y)$ along the line $x=2ay$, which is inside this region.
Along this line, the sign of the vertex condition is positive, at least for $y$ small, as the Taylor expansion 
of the vertex condition is: 
$V_h(2ay,y)=1152(a^7+a^9)y^4+O(y^5)$ and $a>0$.

We can complete the sign of $\kappa'$ in all other regions by just alternating it before and after a vertex branch. 
This completely describes the growth of $\kappa$ on the level sets of $f$ in the plane.
We shall always orient the branches of $f=k$ as indicated in Figure~\ref{fig:growthkappa}, and
in this orientation $\kappa$ will have a definite maximum or a minimum at a
 given vertex. 
\begin{figure}[h!]
  \begin{center}
  \leavevmode
\epsfysize=2in
\epsffile{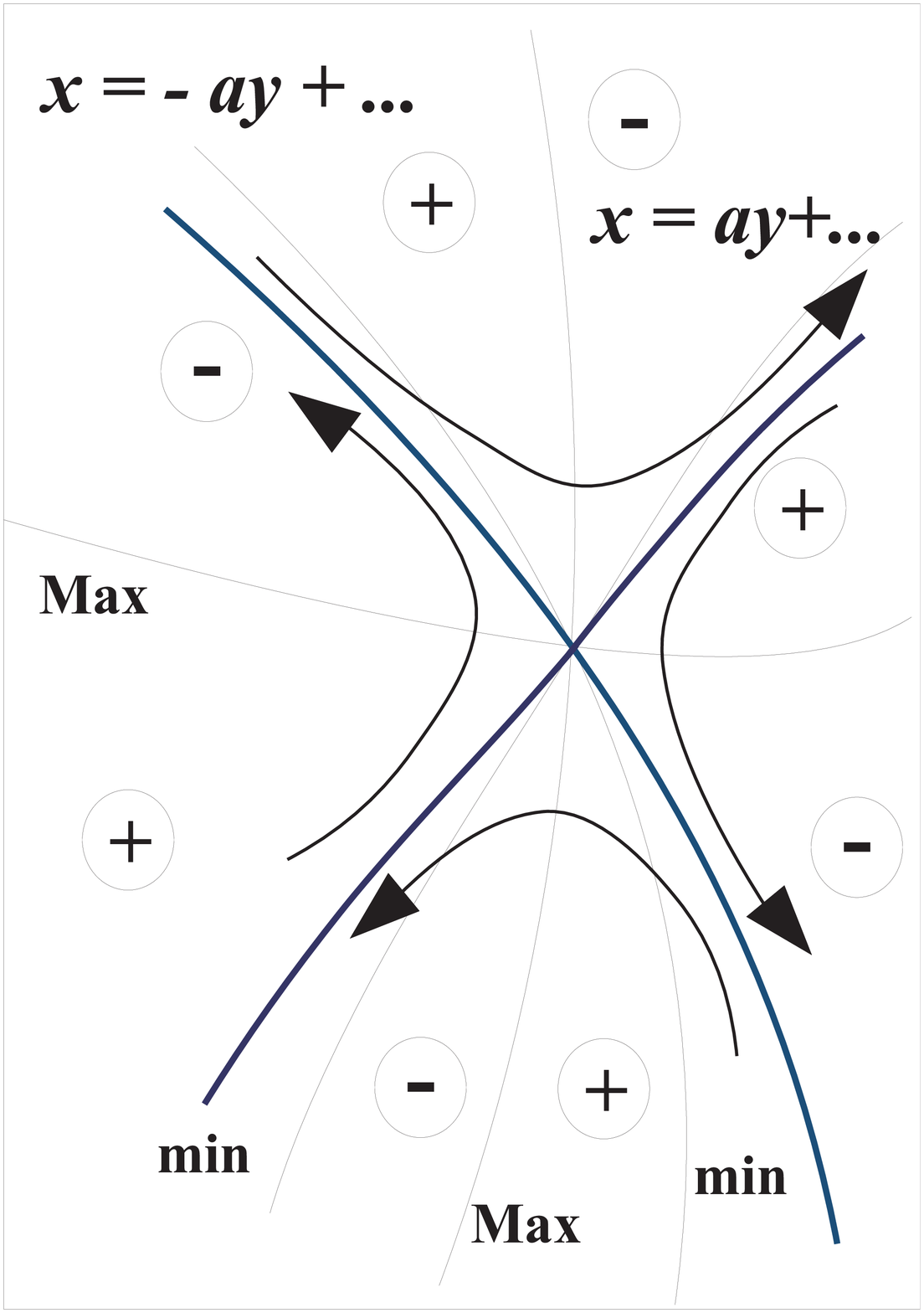}
\epsfysize=2in
\epsffile{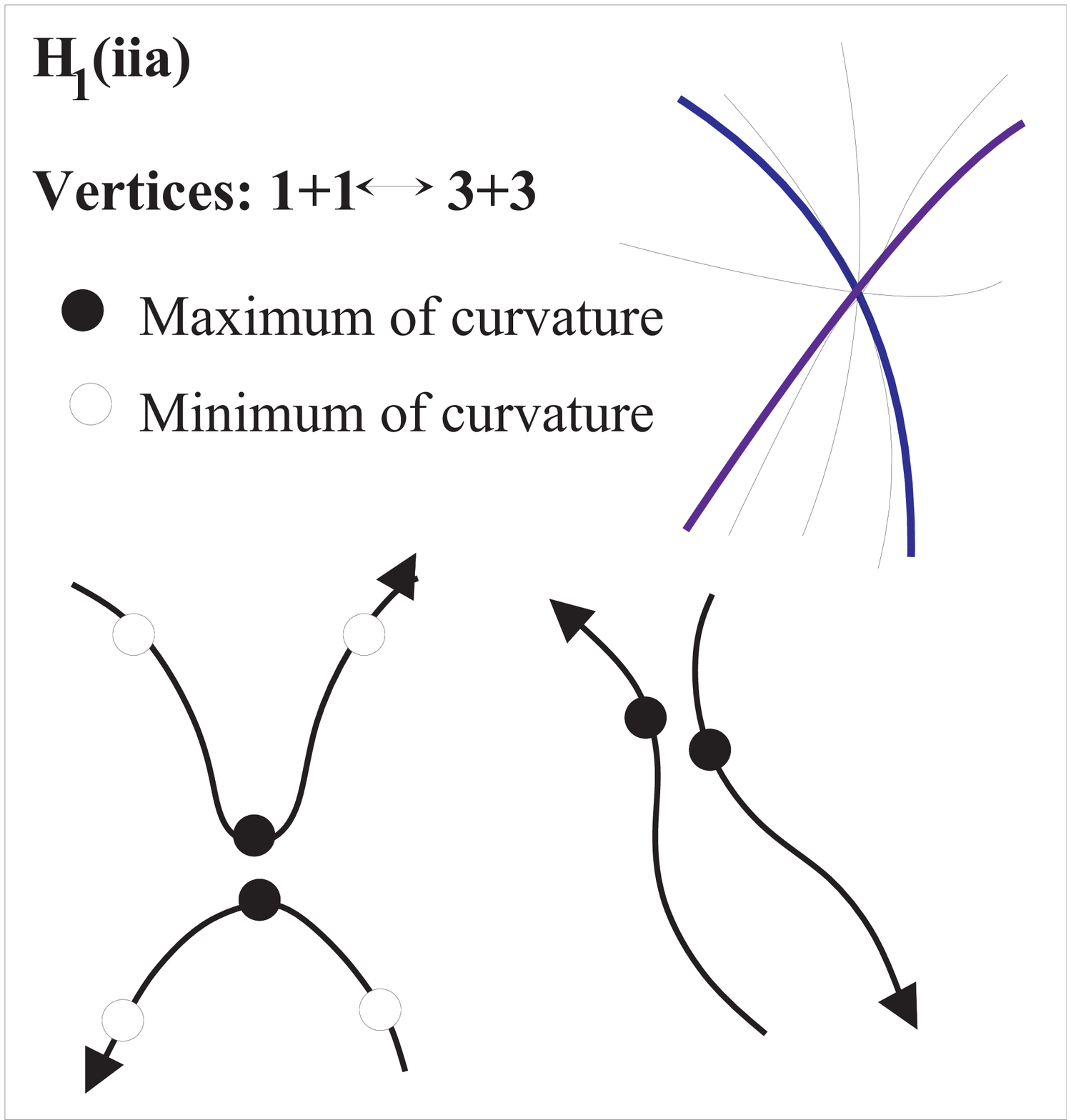}
\end{center}
  \caption{ \small Case ${\bf H}_1$(iia) (compare Figure~\ref{fig:hyp-geom-vertex-inflex}). The sign of $\kappa'$: 
following the indicated orientations on $f=k$, before the 
curve $f=k$ intersects a vertex branch `Max', the derivative $\kappa'$ of $\kappa$ is positive, then vanishes at the 
vertex branch and becomes negative afterwards. So the intersections of the vertex branches ` Max' and the 
curves $f=k$ are the vertices on $f=k$ where $\kappa$ reaches a local maximum; likewise the `min' describes 
the patterns of the local minima of curvature of $f=k$ when $k$ goes through zero. The diagram on the right takes into
account inflexions on $f=k$.}
  \label{fig:growthkappa}
  \end{figure}

\begin{proposition} In the notation of Proposition~\ref{prop:vert-inf},
the limiting curvature of the level curves  $f=k, \ k\to 0$ at vertices  on the various branches is, up to sign, \\
$\bullet$ \ infinite, along $VH_1$ and $VH_{2}$ \\
$\bullet$ \ 
    $f^{(3)}(a)/a(1+a^2)^{3/2}$,  along $VH_{3}$ \\
$\bullet$ \
   $-f^{(3)}(-a)/a(1+a^2)^{3/2}$, along $VH_{4}$.
\label{prop:curvatures}
\end{proposition}

To prove this, we use the Taylor expansions of the branches of the vertex set, given
above in Proposition~\ref{prop:vert-taylor}, and the formula (\ref{eq:kappa2}) for the square
of the curvature of a plane curve.  For the branch $VH_1$ we find the numerator and denominator
of $\kappa^2$ come to $64a^8y^4+O(y^5)$  and $64a^{12}y^6+O(y^7) $ respectively, so
that as $y\to 0$ the limiting curvature is infinite. The situation for $VH_{2}$ is similar.

For $VH_{3}$ the numerator and denominator come to
$64a^4f^{(3)}(a)^2y^6+O(y^7)$ and $64a^6(1+a^2)^3 y^6+O(y^7)$,
which gives the required result. Note that this limiting curvature is zero precisely for a {\em flecnodal}
point, at which the quadratic and cubic terms have a common factor $x-ay$.
The limiting curvatures for both branches $VH_{3}$ and $VH_{4}$ are zero when
the whole of the quadratic terms are a factor of the cubic terms, that is for
the intersection of two flecnodal curves corresponding to different asymptotic
directions on the surface $M$.

\subsection{The vertex transition (VT) set}\label{ss:VT}

Given a generic surface $M$, we can apply our analysis to any point {\bf p}
of the surface: we are then looking at the family of plane sections of the surface close
to the tangent plane section.
The `4-point contact condition' (\ref{eq:VT1}) or (\ref{eq:VT2})
 is generically expected to hold for points
{\bf p} along a set of {\em curves} on $M$, the vertex transition (VT) set.  Of course the VT set
lies entirely in the hyperbolic region, though it may have limit points on the parabolic set (see below);
it separates those points where the family of sections parallel to the tangent plane exhibits behaviour ${\bf H}_1$(i)
in Table~\ref{table:hyper2} from those exhibiting ${\bf H}_1$(ii).

It is clearly of interest to determine, for a given surface $M$, the subregions
into which the hyperbolic region is separated by the VT set.
This set can self-intersect, when {\em both} local branches of 
 $f=0$ have 4-point contact with the corresponding local branches of
the vertex set $V_h=0$: this is $\H_3$ in Table~\ref{table:hyper1}. Also there are 
special points on the VT set where a branch of the vertex set and $f=0$ have
5-point contact: this is $\H_4$ in Table~\ref{table:hyper1}. Although the local conditions are quite easy to calculate---see the
above formulae---it is not so easy to take a global surface and determine the VT set.
We consider below the case of a surface of revolution $M$, which turns out to be non-generic in the sense that
a point of $M$ lies on both branches of the VT set or on neither. We also consider the limit points of the VT set on the parabolic
curve of a general surface $M$.


\medskip\noindent
{\bf Torus and surface of revolution} \ Consider a torus of revolution $M$ in $\R^3$, obtained by
rotating a circle about a line in its plane, not intersecting the circle. Naturally the
VT set will be one or more circular `latitude parallels' of the torus
in view of the circular symmetry. In fact, for a circle of radius $r$ rotating so that its
centre describes a circle $C$ of radius $R>r$, the two latitude
parallels in the hyperbolic region of the torus making an angle $\cos^{-1}(r/R)$ with the plane of $C$ 
lie on both branches of the VT set. Thus for points {\bf p} on these two latitude parallels, {\em both} branches
of the local intersection of $M$ with its tangent plane have 4-point contact with
the corresponding branches of the vertex set at {\bf p} (or equivalently both branches of
the intersection of $M$ with its tangent plane have a vertex at {\bf p}). 
At other hyperbolic points of $M$ {\em neither}
branch has these properties. Crossing the VT set we therefore cross it twice, so that, apart from points {\bf p}  on the VT set itself, 
the pattern of vertices on sections of $M$ parallel to the tangent plane at {\bf p} is always the same.  
In fact we find that, in the expansion of the torus in Monge form at any hyperbolic point, the coefficients
$b_1, b_3, c_1, c_3$ are all zero. It is clear that, in this situation,
the two expressions in Proposition~\ref{prop:VT} become identical so that, in the case $\H_1$ of
Theorem~\ref{vertices-inflexions}, only (ii) is possible. Thus all hyperbolic points away from the VT set
exhibit the same pattern of vertices. Interestingly, when we consider inflexions, then both possibilities
in Table~\ref{table:hyper2} occur. In fact let {\bf p} be a point of the torus of the form
$(r\sin t, 0, R+r\cos t)$ (where the axis of rotation is the $x$-axis and we can without loss of generality
take {\bf p} to be in the $xz$-plane). Then using Proposition~\ref{prop:H1} we find that if
$-r/R<\cos t<0$ then the inflexion transition is $1+1\lr 2+0$ but if $-1<\cos t <-r/R$ it is $2+0\lr 1+1$.
Note that $\cos t<0$ since {\bf p} is hyperbolic, and $t=\pi$ gives the symmetrical case $\H_7$(ii) of Table~\ref{table:hyper2}
and Figure~\ref{fig:torus}.

The same happens in fact for any surface of revolution generated by rotating a plane curve,
say in the $x,z$-plane, about the $z$-axis. We find that $b_1, b_3, c_1, c_3$ are all zero and
the conclusion follows as before. If we rotate the curve $y=0, \ x=a+bz+cz^2+dz^3+ez^4+\cdots$
about the $z$-axis then the condition for the point $(a,0,0)$ to be hyperbolic is $ac>0$ and
the condition for this point to lie on the VT set 
determines $e$ uniquely in terms of $a, b, c, d$.  For example, the curve
$x=a + cz^2 -(c^2/2a)z^4 $ has the latter property, as does
$x=4-2z+2z^2+z^3$.

\medskip\noindent
{\bf The VT set and the parabolic curve} \ The analysis of sections parallel to the tangent plane
at a parabolic point and at a cusp of Gauss is given in \S\S\ref{s:parabolic},~\ref{s:cog}.
Here we are concerned with the hyperbolic region near a parabolic point and we ask which
type, ${\bf H}_1$(i) or ${\bf H}_1$(ii), the points of this region can be.

Suppose we consider a sequence of hyperbolic points
tending to a parabolic point {\bf p} of $M$. If we let $a\to 0$ in (\ref{eq:VT1}) and (\ref{eq:VT2}), the left-hand
sides both tend to $-b_3^2$, since $f^{(3)}(0)=b_3$. Hence, if $b_3\ne 0$ at a point {\bf p} of the parabolic curve,
then all hyperbolic points sufficiently close to {\bf p} are of type ${\bf H}_1$(iia), by Proposition~\ref{prop:H1}(ii).
In particular the VT set cannot have a limit point on the parabolic curve except where $b_3=0$, that is at
the cusps of Gauss.  It is possible to calculate the local form of the VT set at cusps of Gauss; we find the following.\\
$\bullet$ \ At an {\em elliptic} cusp of Gauss {\bf p} (in (\ref{eq:monge}) ${\bf p}=(0,0,0)$ and $\kappa_2=0$, $b_2^2<2\kappa_1c_4$
or, scaling $\kappa_1$ to 2, $b_2^2<4c_4$), there is locally no VT set. \\
$\bullet$ \ At a {\em hyperbolic} cusp of Gauss {\bf p} (the previous inequalities are reversed), there is either  locally no
VT set, or  locally a VT set consisting of two curves tangent to the parabolic curve at {\bf p} and having inflexional
contact with each other (equivalent by a change of coordinates in the parameter plane of $M$ to $(x-y^3)(x+y^3)=0$).
The criterion separating these cases is the sign of a polynomial in coefficients of the Monge form of $M$ at {\bf p}
of order $\le 4$, together with $d_5$. When $d_5=0$ a VT set exists if and only if $c_4$ lies between 0 and
$20b_2^2c_3(b_1b_2-c_3)/(4b_1b_2+c_3)^2$. There is a similar, slightly more complicated formula, for
general $d_5$.


\section{Elliptic points}
\label{s:elliptic}
We sketch this case for completeness; the chief interest for us lies in the symmetry set and medial axis in
the umbilic case as in~\cite{scalespace05}.

At an elliptic point, say $\pp =(0,0,0)$ on a surface $z=f(x,y)$, the two principal curvatures are of the same sign, say
 positive: $\kappa_1>0,\kappa_2>0$. Using (\ref{eq:monge}), the function $f$,
after scaling of the variables $x,y,z$, is of the form
 $f_e(x,y)=x^2+a^2y^2 + b_0x^3+b_1x^2y+b_2xy^2+b_3y^3+$ h.o.t.,
where we may assume $a>0$. We can distinguish two cases here: the generic case where $a\neq 1$ and 
the case $a=1$ of umbilic points, where the principal curvatures are equal. Umbilic points are isolated points in 
the elliptic region of a  surface.  See Figure~\ref{fig:umbilic-geom-vertex-inflex} for the vertex set and
some level curves $f_e(x,y)=k$ in the umbilic case.

\subsection{Proof of Theorem \ref{vertices-inflexions}: elliptic case}\label{ss:proof-elliptic-vert-inflex}
The results on vertices in this case are well-known; to deduce them from the function $V_e$ note that it
has 4-jet
\[-192a^4xy(a^2-1)(x^2+a^2y^2)\]
so that, when $a\ne 1$, there can be only two real branches of $V_e=0$, with tangents $x=0$ and $y=0$.

The inflexion condition $I_e$ has 2-jet $8a^2(x^2+a^2y^2)$, hence the set
$I_e=0$ contains no real points apart from the origin.

The sections $f_e=k$ will therefore have four vertices for small $k>0$, just as in the case of an ellipse.

\smallskip\noindent
{\bf The umbilic case.}
Let us consider the case $a=1$.
The vertex set $V_u=0$ is now given by
\begin{eqnarray*} 
\frac{1}{192}V_u&=&px^5-3qx^4y-2px^3y^2-2qx^2y^3-3pxy^4+qy^5 + \ \mbox{h.o.t.}, \\
&=& (x^2+y^2)(px^3-3qx^2y-3pxy^2+qy^3)+ \ \mbox{h.o.t.},
\end{eqnarray*}
where $p=b_3-b_1, \ q=b_2-b_0$.  The discriminant of the form of degree 3 is
\[ 108(p^2+q^2)^2,\]
so that, unless $p=q=0$ (which amounts to saying that $x^2+y^2$ is a factor of the cubic terms), the
discriminant is $>0$ and the branches of $V_u=0$ through the origin are distinct
and exactly three of them are real. It follows that there are always six vertices on the
section $f_u=k$ for small $k>0$. (Compare~\cite[\S15.3]{porteous}.)  Not surprisingly, there are no inflexions on $f_u=k$. The inflexion
set has equation $I_u=0$, which has the form $8(x^2+y^2)+$h.o.t.

Naturally, the curvature at the vertices tends to infinity as $k\to 0$ through positive values; in fact the curvature behaves
like that of a circle of radius $\sqrt{k}$.

\begin{figure}[h!]
  \begin{center}
  \leavevmode
\epsfysize=1.5in
\epsffile{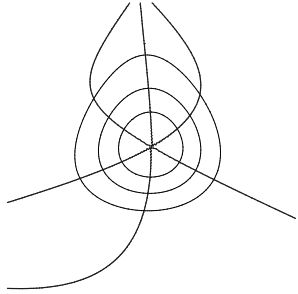}
\end{center}
  \caption{\small Loci of vertices  in a 1-parameter family of  level sets $f=k$ (closed curves), in the umbilic case. The vertex curve
has three branches through the origin, giving rise to six vertices on the level set for all small $k$.}
  \label{fig:umbilic-geom-vertex-inflex}
  \end{figure}

\section{Parabolic case}\label{s:parabolic}
At a parabolic point \pp\ the contact of the surface $M$ with its tangent plane
 is of type $A_2$ at least; we consider the case of ordinary
parabolic points where the contact is exactly $A_2$ in this section.  
One of the principal curvatures vanishes. After scaling of the variables
 $x,y,z$ in (\ref{eq:monge}) $f$ can be written
\begin{equation}
f_p(x,y) = x^2+b_0x^3+b_1x^2y+b_2xy^2+b_3y^3 +\mbox{higher
order terms}
\label{eq:para}
\end{equation}
where $b_3\ne 0$. (The case $b_3=0$ is that of a cusp of Gauss; see \S\ref{s:cog}.)

\begin{proposition} \label{prop:para-vert-inflex}
(i) The vertex set $V_{p}=0$ has three branches, one being smooth and the other two having ordinary cusps. \\
(ii) The inflexion set $I_p=0$ has two branches, one smooth and one having an ordinary cusp.\\
(iii) The zero level set $f_p=0$ has one branch, having an ordinary cusp.
\end{proposition}
See \S\ref{ss:proof-para-vert-inflex} for the proof.

By the same method as in Proposition~\ref{prop:hyp-inflexions} we can show the following.
\begin{proposition}\label{prop:para-vertices}
Suppose that $b_3>0$ (see the Remark below for the contrary case $b_3<0$).\\
{\rm (i)} \ The smooth branch $VP_1$ of the vertex set has the following 3-jet:
 \newline\noindent  
$\bullet$ $VP_1$:  $( -\frac{1}{2}b_2t^2+\frac{b_2(b_1-3b_3)-c_3}{2}t^3,t)$. 

 \smallskip\noindent  
  The two cusped branches of the vertex set have the following 4-jets:
\newline\noindent 
$\bullet$ $VP_{2}$:   $( x_3't^3-\frac{1}{2}b_2t^4 ,-t^2)$,
 \newline\noindent 
$\bullet$ $VP_{3}$: $( x_3''t^3-\frac{1}{2}b_2t^4,-t^2)$,
 \newline\noindent 
 where
 $2x_3'=\sqrt{9+3\sqrt{3} }\sqrt{b_3} $,  and $2x_3''=\sqrt{9-3\sqrt{3}} \sqrt{b_3}$.\\
{\rm (ii)} \ The branches of the inflexion set can be parametrized as
\newline\noindent
 $\bullet$  $(3b_3t, -b_2t+\cdots)$ (recall $b_3\ne 0$)
 \newline\noindent 
 $\bullet$ $(\frac{1}{2}\sqrt{3b_3}t^3 -\frac{3}{8}b_2t^4+ \cdots, -t^2)$.\\
{\rm (iii)} \ The level set $f_p=0$ has the following 5-jet:\\
$\bullet$ $( \sqrt{b_3}t^3-\frac{1}{2}b_2t^4 + \frac{b_2^2+4b_1b_3-4c_4}{8\sqrt{b_3}}t^5$, $-t^2)$.
\end{proposition}
Comparing the coefficients of the $t^3$-terms of the cuspidal branches in (i), (ii) and 
(iii) we have $\frac{1}{2}\sqrt{3b_3}<x_3'' < \sqrt{b_3}
 < x_3'$.  It follows that the branch of $f_p=0$  is always between the 
two cusped branches of  the vertex set, and also the cusped branch of the inflexion set is
inside all these three cusps. See 
Figure~\ref{fig:infl-vert-f-parabolic}.

Hence each level curve $f_p=k$ has only three vertices, near the origin, for  small $k$. Thus, 
when $k$ passes though $0$, the number of  vertices of the curves $f_p=k$ remains unchanged: 
 $3 \leftrightarrow 3$, as claimed in Theorem~\ref{vertices-inflexions}.
The number of inflexions does not change as $k$ passes through $0$: each curve $f=k$ has two inflexions 
near the origin.  Hence the transition of inflexions is $2\leftrightarrow 2$.

\smallskip\noindent
{\bf Remark} \ 
 If $b_3<0$, then in Proposition \ref{prop:para-vertices}  we use $y=t^2$ instead of $y=-t^2$ and replace
$\sqrt{b_3}$ by $\sqrt{-b_3}$ wherever it occurs. The two cases
$'$ and $''$ in (i) are then reversed.

\begin{figure}
\begin{center}
\leavevmode
\epsfysize=1.5in
\epsffile{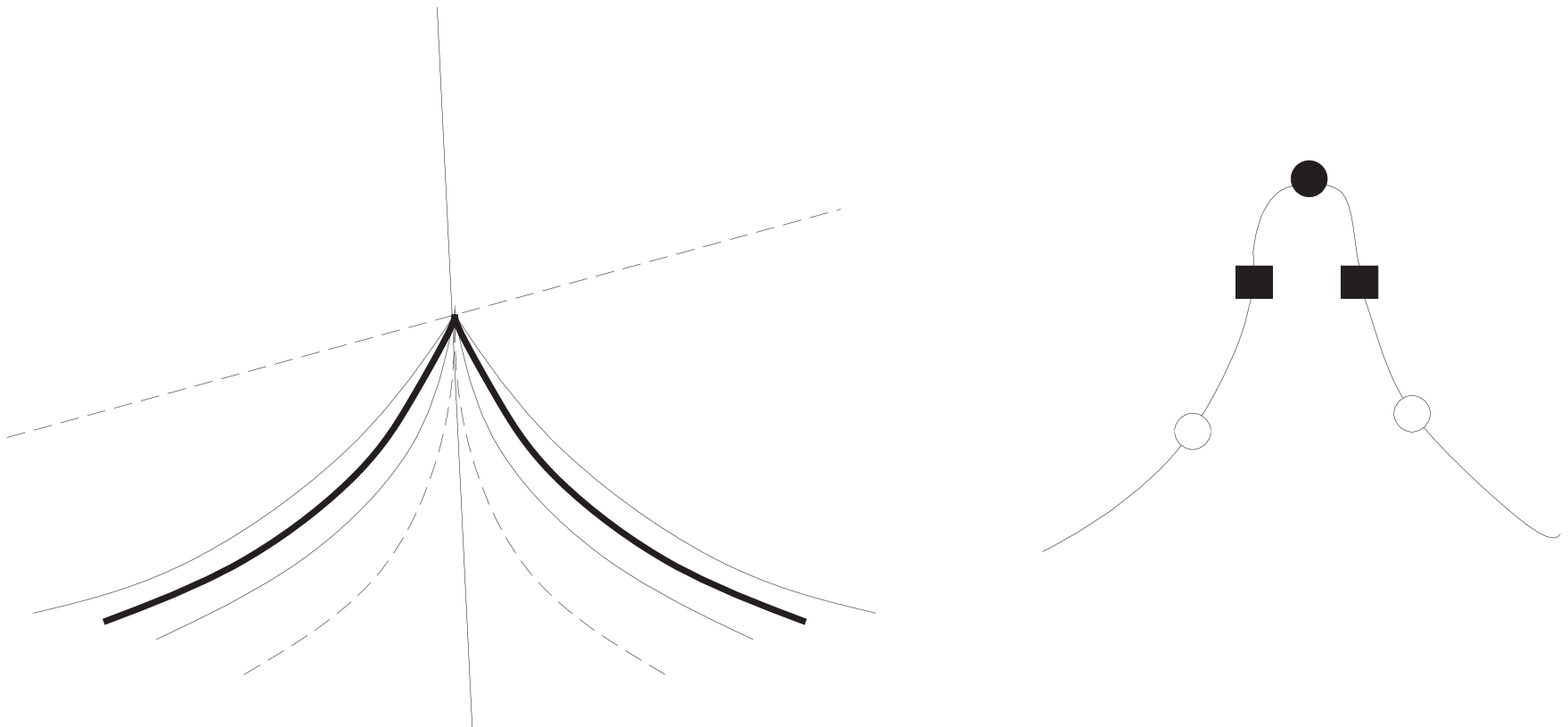}
\hspace*{0.5in}
\epsfysize=1.5in
\epsffile{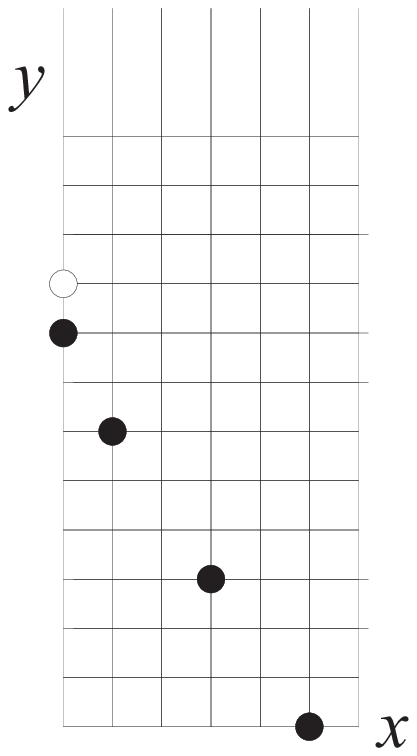}
\end{center}
\caption{\small Left: a schematic picture of the vertex set $V_p=0$ (thin solid line), the inflexion set $I_p=0$ (dashed line) and
the zero level set $f_p=0$ (thick line) in the parabolic case. The vertex set has two cuspidal branches and one
smooth branch, the inflexion set has one cuspidal branch and one smooth branch, and $f_p=0$ has one
cuspidal branch. The level set $f_p=k$ then evolves so that the number of vertices remains as 3 and the number of inflexions
as 2 for both signs of $k$, with $k$ small. Centre: a sketch of the level curve $f_p=k$ for $k\ne 0$, marking
vertices (circles) and inflexions (squares). Right: the Newton polygon for the parabolic case; see \S\ref{ss:proof-para-vert-inflex}.}
\label{fig:infl-vert-f-parabolic}
\end{figure}

\subsection{Proof of Proposition~\ref{prop:para-vert-inflex}}\label{ss:proof-para-vert-inflex}
 In this section we show briefly how we are able to deduce that the vertex and inflexion sets have
 branches as claimed above. We  do so by looking at the Newton polygon and then applying the well-known
techniques of blowing-up combined with the implicit function theorem.  
We give this example in detail; all the other cases encountered in this
article can be dealt with similarly.

The Newton polygon for the function $V_p$ contains 
the following monomials with coefficients:
$192b_3x^5+864b_3^2x^3y^3+648b_3^3xy^6+324b_2b_3^3y^8.$
Since $b_3\ne 0$ all but
the last term are definitely present. The last term is absent when $b_2=0$, which means that the
parabolic curve is tangent to the {\em other} principal direction at the origin. When this is the case, there is a
term $324b_3^3c_3y^9$, which will be present unless $c_3=0$. Generically $b_2=c_3=0$ will not happen anywhere on our
surface $M$. Thus the Newton polygon has terms $x^5, \ x^3y^3,\ xy^6 $ and either $y^8$ or $y^9$;
see Figure~\ref{fig:infl-vert-f-parabolic}. 

Let us write the above as $g(x,y)=ax^5+bx^3y^3+cxy^6+dy^8$ and consider the case where $d\ne 0$.  Note that
$a, b, c$ and $b^2-4ac=248832b_3^4$ are all nonzero. The function $V_p$ will then be of the form $g \ +$ terms above the Newton polygon;
we can think of the latter as linear combinations of monomials $x^my^n$ where $3m+2n>15$ and $(m,n)\ne (0,8)$.
We first blow up by $x=ty$, so that the `blow-down' transformation is $(t,y)\to (ty,y)$ and $y=0$ is the exceptional
divisor. The result after cancelling $y^5$ is
\[ at^5 + bt^3y + cty^2 + dy^3 + \ \mbox{linear combination of monomials} \ t^my^{m+n-5}.\]
Note that $m+n>5$ for all monomials above the Newton polygon. Hence intersecting with $y=0$ gives five
coincident points at the origin $t=y=0$ (and there are no points sent to infinity, that is $y=tx$ produces no
points on the exceptional divisor, using $a\ne 0$).

For the second blow-up we use  $y=ut$, with 
blow-down map $(t,u)\to (t,ut)$; we find after cancelling $t^3$
\[ at^2+but+cu^2+du^3 + \ \mbox{linear combination of monomials} \ u^{m+n-5}t^{2m+n-8},\]
and $2m + n-8>0$ for all monomials above the Newton polygon. This meets $t=0$ in 
$cu^2+du^3=0$, that is a cusp at the origin and a transverse crossing of the $u$-axis at $u=-c/d$, since
$d\ne 0$.  The transverse crossing provides us with a smooth branch of the blown-up curve $V_p=0$,
parametrized by $t$, using the implicit function theorem, and by blowing-down we obtain one of
the branches of our curve $V_p=0$ (in fact also a smooth branch).  No further points are obtained
from the alternative blow-up $t=uy$.

Blowing up the origin a third time, using $t=uw$, with blow-down map $(w,u)\to (uw, u)$, we obtain after
cancelling $u^2$,
\[ aw^2+bw+c+du + \ \mbox{linear combination of monomials} \ u^{3m+2n-15}w^{2m+n-8},\]
and again $3m+2n-15>0$ for all points above the Newton polygon. Finally this meets the exceptional
divisor $u=0$ in distinct points, since $b^2\ne 4ac$, each of which gives a transverse crossing of
$u=0$ so that the two branches of the blown-up curve can be locally parametrized smoothly by $u$, using
the implicit function theorem. These
branches blow-down to the remaining two branches (actually ordinary cusps) of $V_p=0$.

\subsection{\bf The limiting curvature of $f_p=k$, at vertices.}
Here again, we would like to evaluate the curvature $\kappa$ of  $f_p=k$ at a vertex. Then we will take 
the limit of $\kappa$, as one approaches the parabolic point {\bf p}$=(0,0,0)$, along that vertex branch.
\begin{proposition} The limiting curvature of the level curves $f=k$ is infinite as $k\to 0$, at vertices 
on any of the branches of the vertex set.
\label{prop:lim-curv-parab}
\end{proposition}
We substitute the parametrizations of the branches of $V_p=0$ given in Proposition~\ref{prop:para-vertices}
into the expression for $\kappa^2$ given in (\ref{eq:kappa2}). The result is (in all cases using $b_3\ne 0$)
for the branch $VP_1$, $\kappa^2\sim t^{-4}$, while for $VP_2$ and $VP_3$, $\kappa^2\sim t^{-2}$.
The result follows.

\section{Non-degenerate cusps of Gauss}
\label{s:cog}

For a non-degenerate cusp of Gauss the Monge form (\ref{eq:monge}) can be written, after scaling the variables, as
\[
   f_{g}=x^2+b_0x^3+b_1x^2y+b_2xy^2+c_0x^4+c_1x^3y+c_2x^2y^2+c_3xy^3+c _4y^4 +\  \mbox{h.o.t}
\]   
where
$b_2^2-4c_4 \ne 0$, that is, the lowest degree terms in the weighted sense,
namely the $x^2, xy^2$ and $y^4$ terms, are non-degenerate.  Since cusps of Gauss are isolated we can assume generically that
other conditions on the coefficients are avoided. By changing the sign of $x$ if necessary we can
assume $b_2>0$. There are two broad cases:\\
{\em Elliptic cusp of Gauss}: $b_2^2 - 4 c_4 < 0$. Then 
the curve $f_g = k$ is locally a closed loop for $k>0$ and empty for $k<0$. \\
{\em Hyperbolic cusp of Gauss}:  $b_2^2 - 4 c_4 > 0$. Then $f_g=k$ has two local branches for
$k\ne 0$ and two tangential branches for $k=0$. \\
Note that the principal direction $x=0$ is tangent to the parabolic curve at a cusp of Gauss. See~\cite{cog} for an extensive discussion
of cusps of Gauss, and~\cite[pp.245,276]{solid-shape} for further geometrical information.

\subsection{Vertices and inflexions on level sets at a cusp of Gauss}\label{ss:vert-inf-cog}
In a neighbourhood of a cusp of Gauss, the vertex condition now reads
\[V_g=
192(c_3-b_1b_2)x^6+192(4c_4-b_2^2)x^5y+ \mbox{h.o.t.} \]
  Since $b_2^2-4c_4 \ne 0$,  there will be a nonzero coefficient of $x^5y$ here.

\begin{proposition} \label{prop:vert-inf-cog} {\rm (i)} \ In the case of an elliptic cusp of Gauss i.e. $b_2^2 - 4 c_4 < 0$  (the closed 
curve intersection), there are two  smooth  real branches of the vertex set $V_g=0$ through the origin, 
one of which is tangent to the axis $x=0$,  and the other one to $ (b_2^2-4c_4)y =(c_3-b_1b_2)x$.

\smallskip\noindent
{\rm (ii)} \ In the case $b_2^2 - 4 c_4 > 0$ (hyperbolic cusp of Gauss), the vertex set has six smooth
real branches $VG_i$ for $i=1,\cdots 6$. All except $VG_6$ are tangent to $x=0$ while  
$VG_6$ is tangent to $ (b_2^2-4c_4)y =(c_3-b_1b_2)x$.

\smallskip\noindent
{\rm (iii)} \  If in addition to $b_2^2 - 4 c_4 > 0$, we have $b_2^2-8c_4>0$, then the inflexion set has three smooth 
branches (see Figures~\ref{fig:cog-vert-2inflex} and~\ref{fig:cog-vert-3inflex}), whereas when $b_2^2-8c_4<0$, there is 
only one smooth branch (see Figure~\ref{fig:cog-vert-1inflex}).
\end{proposition} 

The claimed number of branches can be deduced from the Newton polygon in the same way as \S\ref{ss:proof-para-vert-inflex};
the present case is easier. The Newton polygon for $V_g$ is illustrated in Figure~\ref{fig:elliptic-cog}, right.
The terms on the Newton polygon are
\begin{eqnarray*}
&&192(c_3-b_1b_2)x^6-192(b_2^2-4c_4)x^5y -480b_2(b_2^2-4c_4)x^4y^3\\
&& -48(b_2^2-4c_4)(7b_2^2+12c_4)x^3y^5 -24b_2(b_2^2-4c_4)(b_2^2+36c_4)x^2y^7\\
&&+24(b_2^2-4c_4)(b_2^4-10b_2^2c_4-16c_4^2)xy^9+24b_2c_4(b_2^2-4c_4)(b_2^2-8c_4)y^{11}.
\end{eqnarray*}
The key fact is this: ignoring the first term and then cancelling $y$, the remaining terms form a
quintic polynomial in $x$ and $y^2$ which has {\em distinct} roots; in fact it factorizes as\\
$(b_2^2-4c_4)(2x+b_2y^2)(x^2+b_2xy^2+c_4y^4)(4x^2+4b_2xy^2-(b_2^2-8c_4)y^4)$. The discriminant 
 is a nonzero constant times $(b_2^2-4c_4)^{18}$ and the number of real roots is
1 for an elliptic cusp and 5 for a hyperbolic cusp. Two blow-ups $x=ty$ and $y=tu$ suffice to
find the real branches of the singular point $V_g=0$.

\medskip
 
We can parametrize the branches of the inflexion set as follows.
Substitute $x=x_1y+x_2y^2+\ldots$ in the inflexion condition $I_g=0$; this gives the solution $x_1=0$,
implying that the inflexion branches are all tangent to the $y$-axis.  Then the coefficient $z$ of $y^2$ is 
a solution of a cubic equation $I(z)=0$ where
\begin{equation}
I(z)=-4c_4(b_2^2-8c_4)-6b_2(b_2^2-8c_4)z+48c_4z^2+8b_2z^3
\label{eq:x2cog-infl}
\end{equation}
with discriminant  
$D=6912 (b_2^2-8c_4)(b_2^2-4c_4)^4$. 
So $D>0$ , giving 3 solutions for $z$, if and only if   $b_2^2-8c_4  >0$. However, if
$b_2^2<4c_4$ (elliptic cusp) then $c_4>0$ so automatically $b_2^2<8c_4$ and the
3 solutions case applies only to hyperbolic cusps.

\begin{figure}
  \begin{center}
  \leavevmode
\epsfysize=1.5in
\epsffile{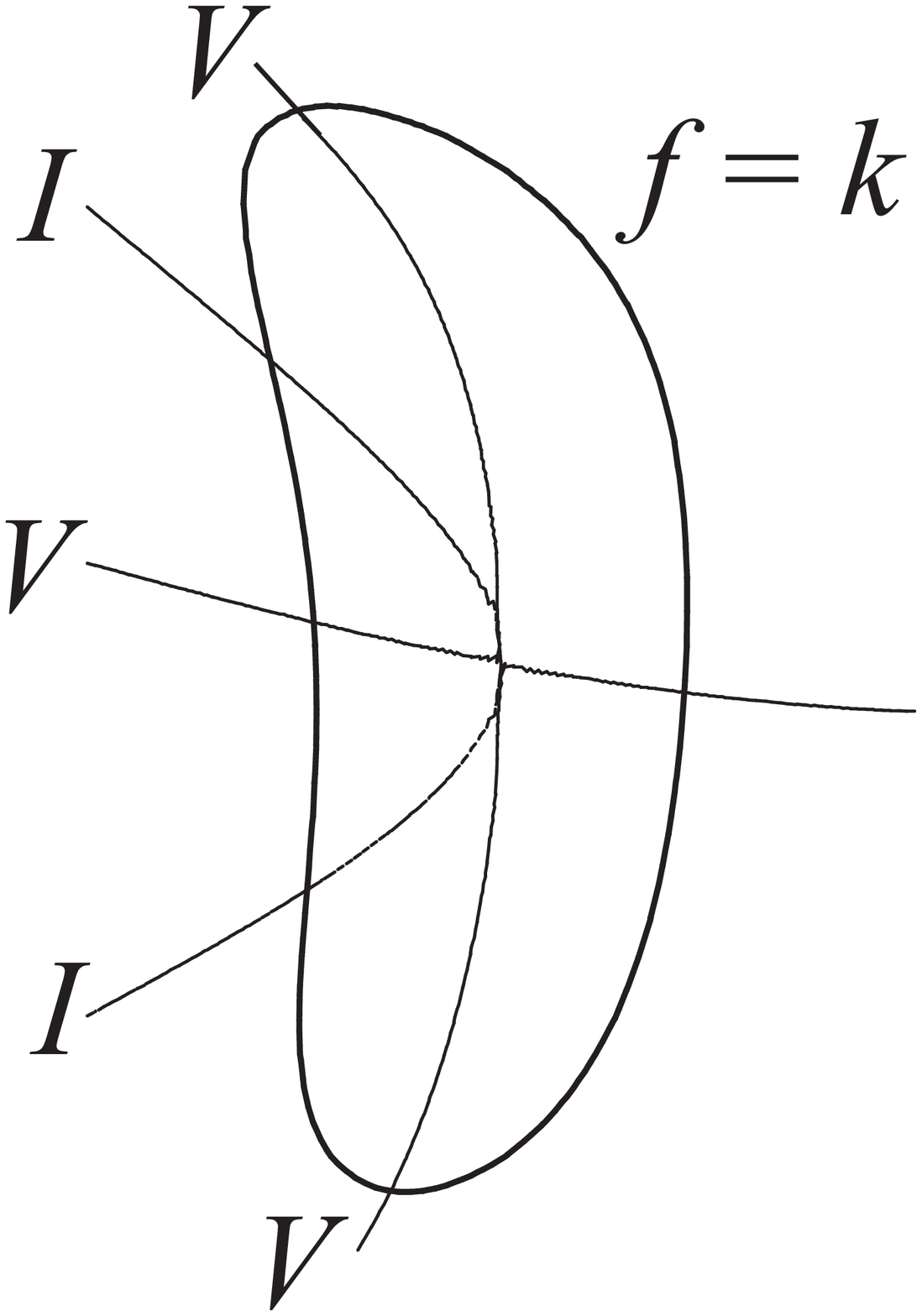}
\hspace*{1in}
\epsfysize=1.5in
\epsffile{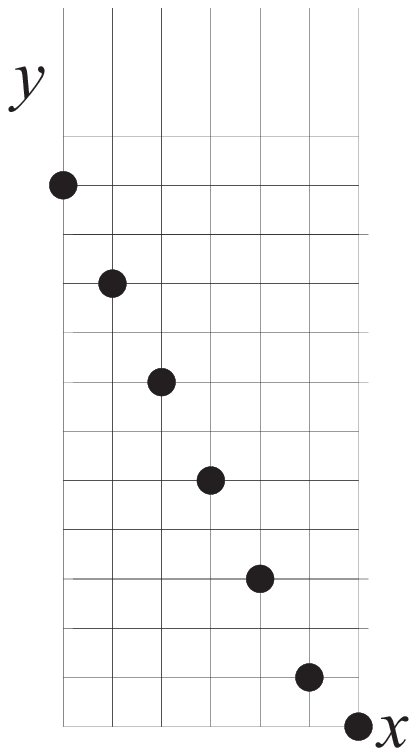}
  \end{center}
\caption{\small Left: Vertices and inflexions in the case of an elliptic cusp of Gauss: 
the curves marked $V$ are the vertex set, those marked $I$ are the inflexion set and
$f=k$ is one level set of $f$. As $k$ increases through $0$, the curve passes from empty
to one with four vertices and two inflexions. Right: the Newton polygon for a
hyperbolic cusp of Gauss; compare \S\ref{ss:vert-inf-cog}.}
\label{fig:elliptic-cog}
  \end{figure}

Once we know the number of smooth branches we can find their Taylor expansions using the same
method of substitution of power series that was used in the previous cases. We find the following.

\begin{proposition} \label{prop:VGFG} When $b_2^2 - 4 c_4 > 0$ (hyperbolic cusp of Gauss), 
 the branches of the vertex set, tangent to the principal direction $x=0$, can be parametrized as follows:\\
$ VG_1$: $x=-\frac{1}{2} \left(b_2-\sqrt{b_2^2-4c_4}\right)y^2 + $ h.o.t., \ \ 
$VG_2$:  $x=-\frac{1}{2} \left(b_2+\sqrt{b_2^2-4c_4}\right)y^2+$ h.o.t.,
\newline\noindent
 $VG_3$: $x=-\frac{1}{2} \left(b_2-\sqrt{2b_2^2-8c_4}\right)y^2 +$ h.o.t., \ \ 
$VG_4$:  $x=-\frac{1}{2} \left(b_2+\sqrt{2b_2^2-8c_4}\right)y^2 +$  h.o.t.\\
$VG_5$: $x=-\frac{1}{2}b_2y^2  -\frac{1}{2}(c_3-b_1b_2)y^3 +$ h.o.t.

\smallskip\noindent
The level set $f=0$ has two branches which can be parametrized as:
\newline\noindent
$FG_1:$ $x=-\frac{1}{2}\left(b_2-\sqrt{b_2^2-4c_4}\right)y^2 +$ h.o.t., \ \ 
 $FG_2:$ $x=-\frac{1}{2}\left(b_2+\sqrt{b_2^2-4c_4}\right)y^2 +$  h.o.t.
\end{proposition}
This Proposition implies in particular that the vertex branch $VG_1$ and the
 branch  $FG_1$ of $f=0$ have at least $3$-point contact at the origin. The same holds for $VG_2$ and $FG_2$. 

The conditions for 4-point contact are given below; since cusps of Gauss are isolated on a generic surface, only
the signs of the expressions below will be of significance.

\begin{proposition} The vertex branch $VG_1$ and the branch  $FG_1$ of $f_g=0$ have at least 
$4$-point contact at the origin if and only if $D_1=0 $ where\\
 $D_1=-b_1b_2^2+b_1b_2\sqrt{b_2^2-4c_4}+2b_1c_4+b_2c_3-c_3\sqrt{b_2^2-4c_4}-2d_5$.

 The same holds for $VG_2$ and $FG_2$ if and only if $D_2=0$ where
\\
$ D_2=b_1b_2^2+b_1b_2\sqrt{b_2^2-4c_4}-2b_1c_4-b_2c_3-c_3\sqrt{b_2^2-4c_4} + 2d_5$.
\label{prop:D1D2}
\end{proposition}
The signs of the $D_i$ determine the relative positions
of the branches $VG_i$ and $FG_i$. More precisely,  $D_1>0$ if and only if,  above the $x$-axis, the curve
$VG_1$ is to the right of $FG_1$. (Below the $x$-axis this is reversed, since they have 3-point contact
at the origin.)  Similarly, $D_2>0$ if and only if, above the $x$-axis, $VG_2$ is to the right of $FG_2$. Note that
both $D_1>0$ and $D_2>0$ can be regarded as conditions on the coefficient $d_5$.

\subsection{Hyperbolic cusp of Gauss}

Let $x_{2i},\  i=1,...,5$ be the coefficient of $y^2$ in the expansion of the branch $VG_i$ as in Proposition~\ref{prop:VGFG},
and let $z_0$ or $z_1 < z_2< z_3$ denote the real roots of (\ref{eq:x2cog-infl}), as appropriate. Thus the $z_i$ are the
coefficients of $y^2$ in the expansion(s) of the branch(es) of the inflexion set: $x=z_iy^2+\ldots$. Recall that
we also assume $b_2>0$. The following is obtained from the expressions in
Proposition~\ref{prop:VGFG} and the sign of the polynomial $I$ in (\ref{eq:x2cog-infl}) at  the values $x_{2i}$.
\begin{proposition} {\rm (a)} Suppose $b_2^2-8c_4  >0$. \\
{\rm  (1)} If  $c_4>0$ then 
$x_{24} < x_{22}  < z_1 < x_{25} < x_{21}<z_2 <0< x_{23}< z_3$.
\\
{\rm (2)} If $c_4 <0$ then $x_{24} < x_{22}  < z_1 < x_{25}<0 <z_2 < x_{21}< x_{23}< z_3$.
\\
{\rm (b)} Suppose  $b_2^2-8c_4  <0$. Then
 $x_{24} < x_{22}  < z_0 < x_{25} < x_{21}<x_{23}<0$.
\end{proposition}

These are illustrated in Figures~\ref{fig:cog-vert-1inflex}-\ref{fig:cog-vert-3inflex}. The sign of the derivative of curvature
is determined as for the hyperbolic case; this determines the pattern
of maxima and minima of curvature. The statements of Theorem~\ref{vertices-inflexions}, case (HCG), follow from
these diagrams.

\begin{figure}
  \begin{center}
  \leavevmode
\epsfysize=1.6in
\epsffile{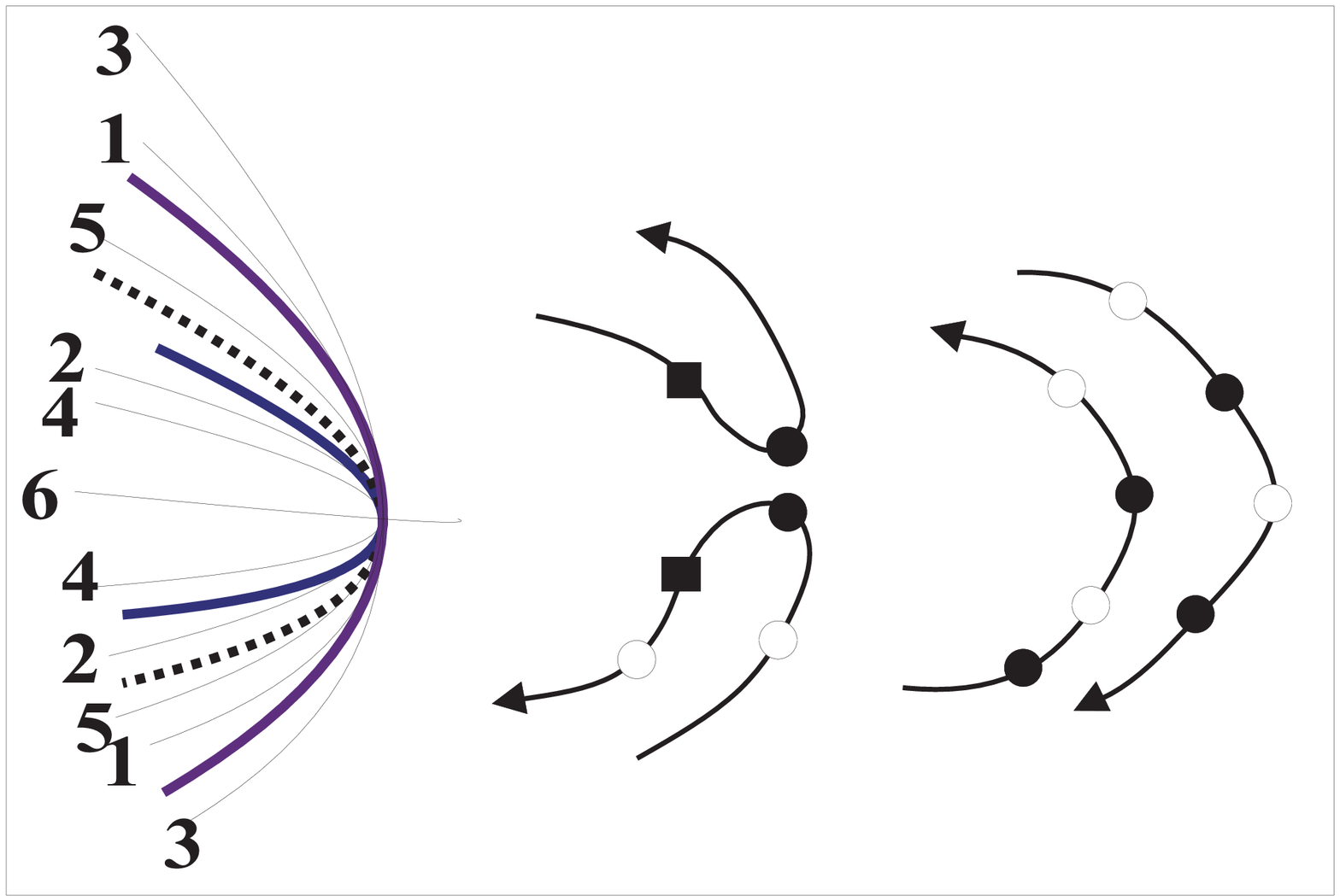}
\epsfysize=1.6in
\epsffile{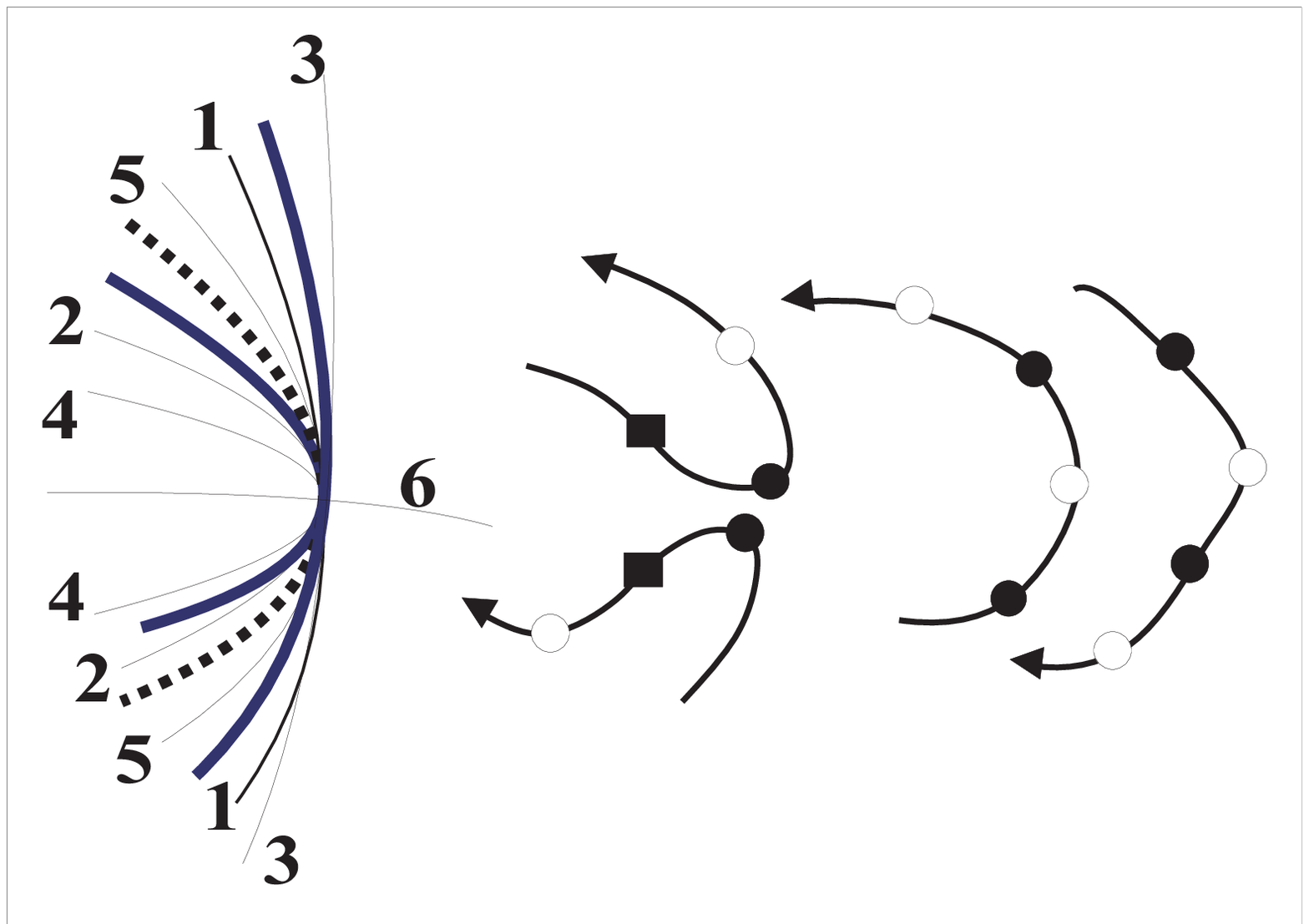}
\end{center}
  \caption{\small Sketch of the hyperbolic cusp of Gauss case, $b_2>0, \ 4c_4<b_2^2<8c_4$ (so $c_4>0$); see Proposition~\ref{prop:vert-inf-cog}.
 The left box has
$D_1$ and $D_2$ as in Proposition~\ref{prop:D1D2} of the same sign (negative) and the right box of opposite signs ($D_1>0$).
The thick lines are $f=0$ on the left of each figure and  $f=k$ for the two signs of small nonzero $k$ on the right. The thin lines are the $VG_i$  of
Proposition~\ref{prop:VGFG}, labelled by $i$, and the dashed line is the single branch of the inflexion
set. As before, solid circles are maxima and open circles are minima of curvature, for the orientations indicated, and squares are inflexions.}
  \label{fig:cog-vert-1inflex}
  \end{figure}

\begin{figure}
  \begin{center}
  \leavevmode
\epsfysize=1.8in
\epsffile{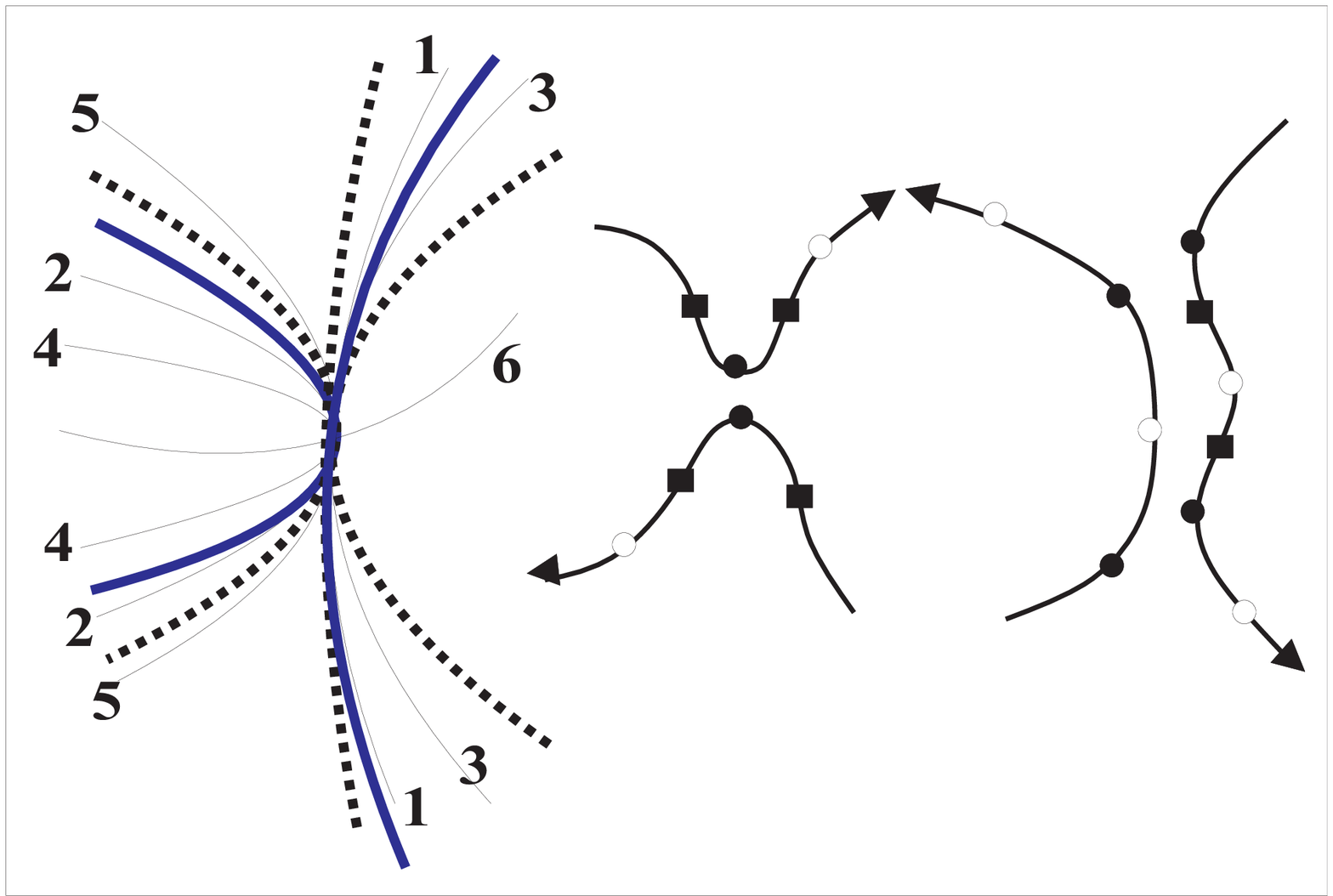}
\epsfysize=1.8in
\epsffile{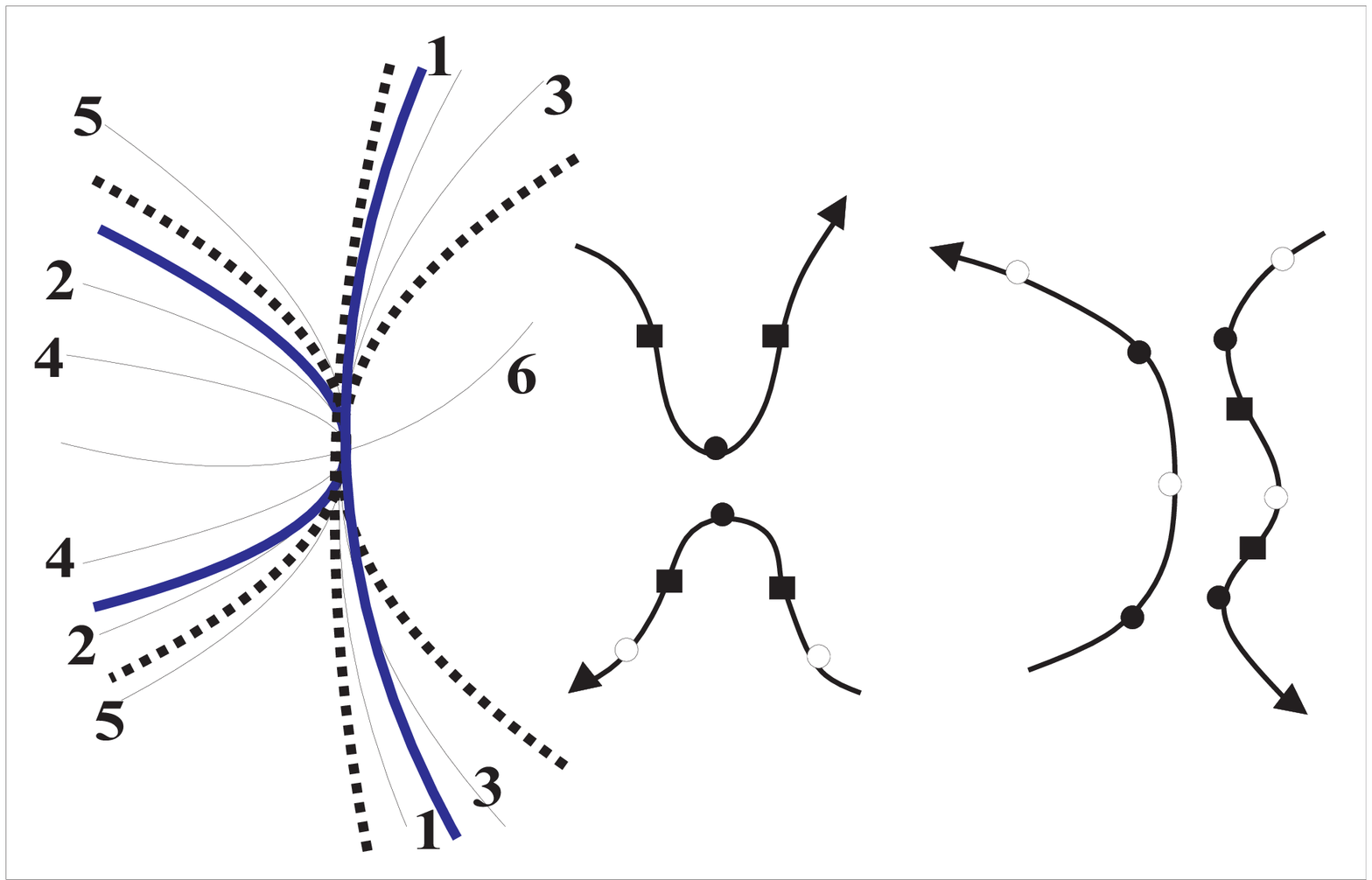}
\end{center}
  \caption{\small Sketch of the hyperbolic cusp of Gauss case, $b_2>0,  \ c_4<0$ (hence $b_2^2>8c_4$); see Proposition~\ref{prop:vert-inf-cog}.
 The left box has
$D_1$ and $D_2$ as in Proposition~\ref{prop:D1D2} of the same sign (negative) and the right box of opposite signs ($D_1>0$).
The thick lines are $f=0$ on the left of each figure and  $f=k$ for the two signs of small nonzero $k$ on the right. The thin lines are the $VG_i$  of
Proposition~\ref{prop:VGFG}, labelled by $i$, and the dashed lines are the three branches of the inflexion
set. As before, solid circles are maxima and open circles are minima of curvature, for the orientations indicated, and squares are inflexions.}
  \label{fig:cog-vert-2inflex}
  \end{figure}

\begin{figure}
  \begin{center}
  \leavevmode
\epsfysize=1.4in
\epsffile{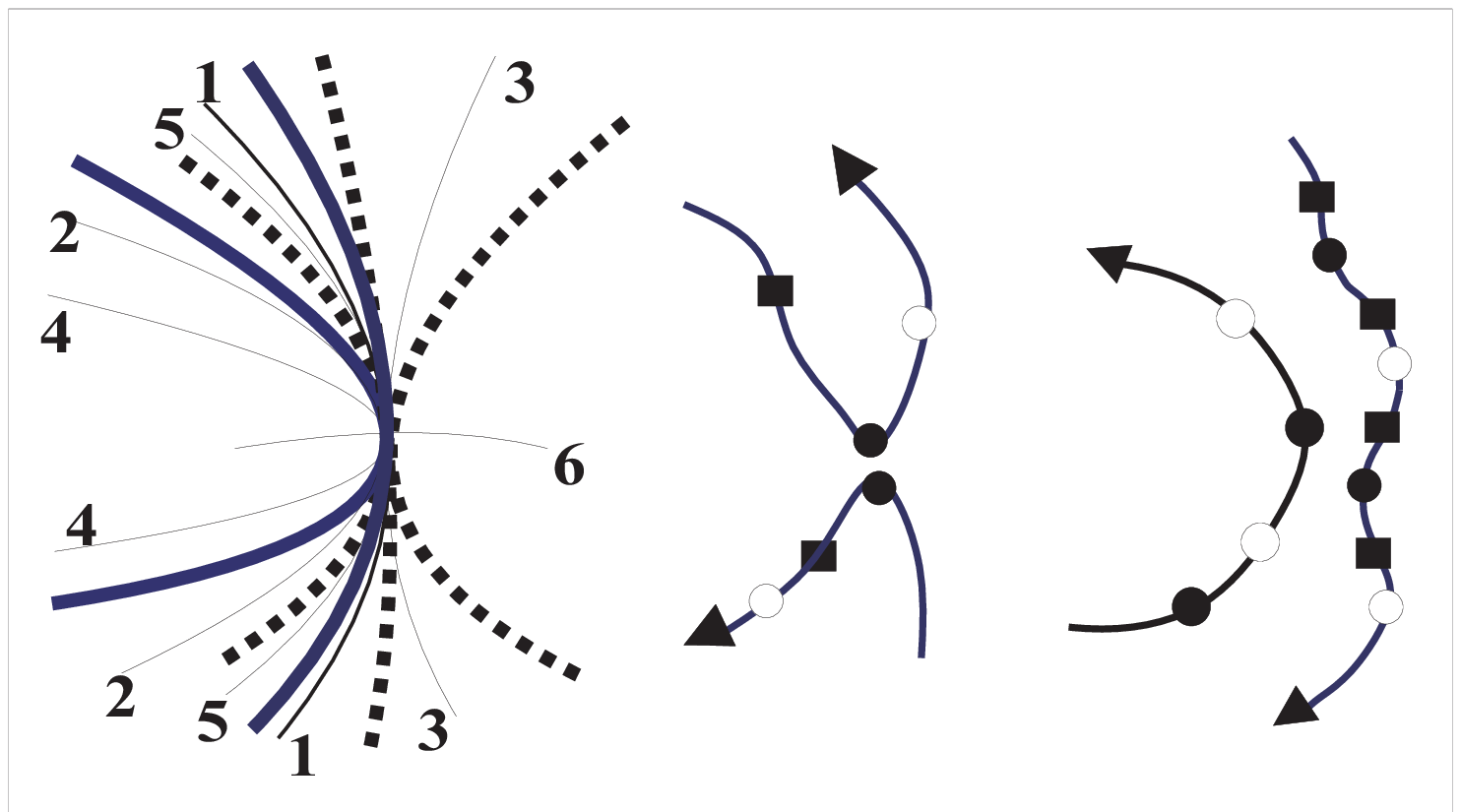}
\epsfysize=1.4in
\epsffile{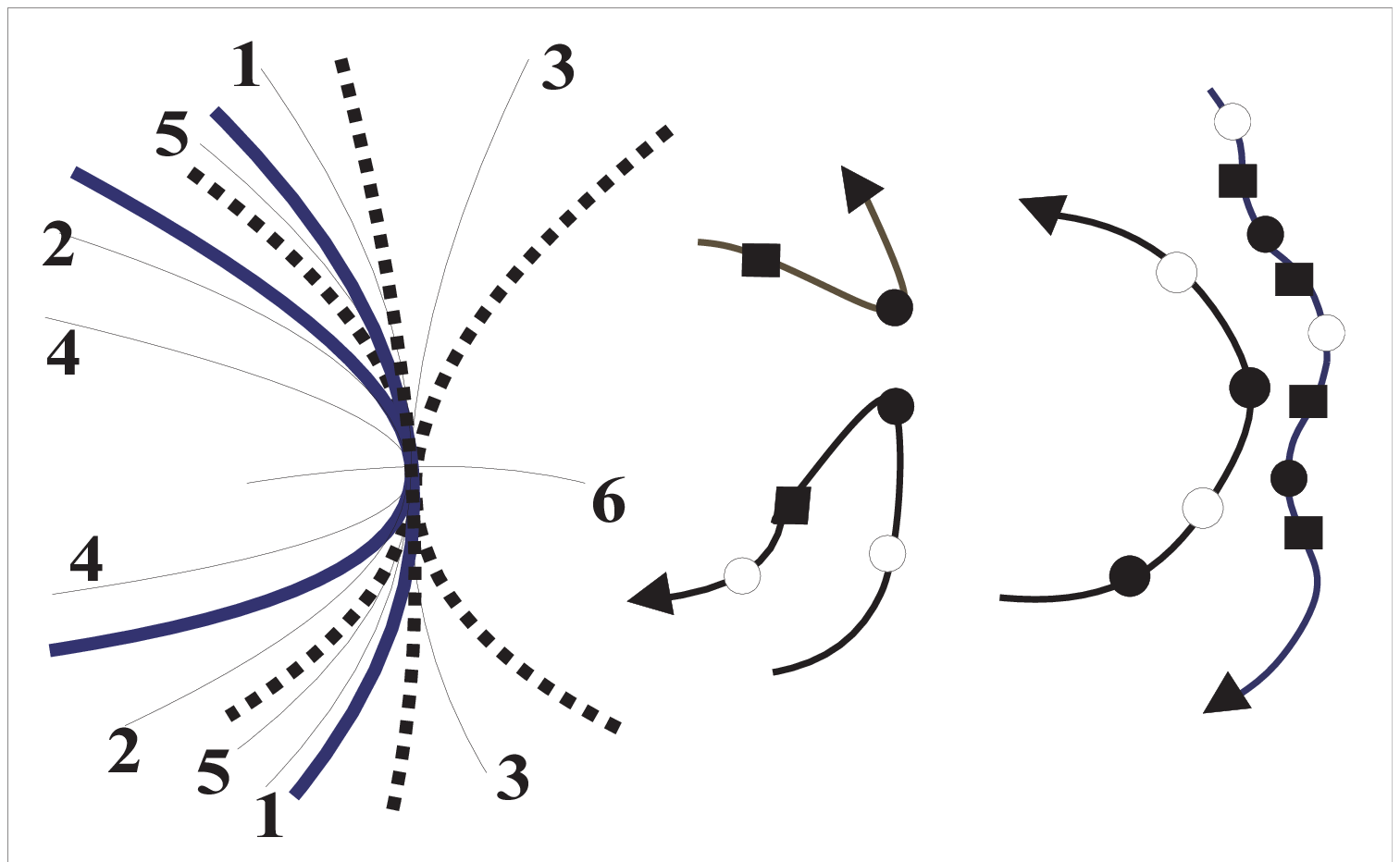}
\end{center}
  \caption{\small As for Figure~\ref{fig:cog-vert-2inflex} except that $b_2>0, \ b_2^2>8c_4, \ c_4>0$; see Proposition~\ref{prop:vert-inf-cog}.}
  \label{fig:cog-vert-3inflex}
  \end{figure}

\section{Conclusion}\label{s:conclusion}
In this article, we have derived detailed results on the pattern of vertices and inflexions on families of plane curves
of the form $f(x,y)=k$, which can be interpreted as the parallel plane sections
of a generic surface close to the tangent plane at a given point \pp. This is part of an investigation of the
symmetry sets and medial axes of 1-parameter families of plane curves which evolve through a singular
member. The symmetry set of a nonsingular plane curve $\gamma$ is the closure of the locus of
centres of circles tangent to $\gamma$ in more than one place (`bitangent circles'). It has endpoints in the cusps of the evolute, that is
at the centres of curvature of the vertices of $\gamma$. Thus the pattern of vertices has a strong influence on the
branches of the symmetry set. Inflexions have a direct effect on the evolute---it goes to infinity---and, through
the associated double tangents, an indirect effect on the symmetry set, which has a point at infinity for every
double tangent (a bitangent circle of infinite radius). The limiting curvatures at vertices, as $k\to 0$, determines
the limiting position of the endpoints of the symmetry set as the plane section becomes singular.

The investigation of symmetry sets involves many other factors, such as an investigation of circles which
are tangent in three places to $\gamma$ (these produce triple crossings on the symmetry set) and circles
which are circles of curvature at one point of $\gamma$ and tangent elsewhere (these produce cusps on
the symmetry set). These and other matters are reported elsewhere, beginning with \cite{scalespace05}.

We conclude with some remarks and questions about the material of this article. \\
1. Is it possible to calculate the VT curve for classes of global examples where the two branches
do not coincide? Compare \S\ref{ss:VT}. \\
2. Can the parabolic and cusp of Gauss cases be approached by more general methods of singularity theory,
as in \cite{garay,uribe-vargas}?\\
3. For the purpose of plotting symmetry sets it is much more convenient to have a parametrized curve rather than
a level set $f(x,y)=k$. A method of parametrizing the level sets to
arbitrarily high accuracy is given in \cite{maths-of-surfaces}.

\end{document}